\documentclass[11pt,leqno,draft]{article}
\usepackage{amsfonts}
\usepackage{latexsym}
\usepackage{amsmath}
\usepackage{amssymb}
\usepackage{amsthm}
\newtheorem{theorem}{Theorem}[section]
\newtheorem{lemma}[theorem]{Lemma}

\newtheorem{proposition}[theorem]{Proposition}
\newtheorem{corollary}[theorem]{Corollary}
\newtheorem{remark}[theorem]{Remark}

\theoremstyle{definition}

\theoremstyle{remark}

\newtheorem*{note*}{Note}

\numberwithin{equation}{section}

\makeatletter
\newcommand{\ls}{\leqslant}
\newcommand{\gr}{\geqslant}
\makeatother

\begin{document}
\small

\title{\bf Geometry of the $L_q$-centroid bodies of an isotropic log-concave measure}

\medskip

\author{Apostolos Giannopoulos, Pantelis Stavrakakis, Antonis Tsolomitis\\[1ex]
 and Beatrice-Helen Vritsiou}

\date{}

\mathversion{bold}
\maketitle
\mathversion{normal}

\begin{abstract}
\footnotesize We study some geometric properties of the $L_q$-centroid
bodies $Z_q(\mu )$ of an isotropic log-concave measure $\mu $ on ${\mathbb R}^n$. For any $2\ls q\ls\sqrt{n}$
and for $\varepsilon \in (\varepsilon_0(q,n),1)$ we determine the inradius of a random $(1-\varepsilon )n$-dimensional projection
of $Z_q(\mu )$ up to a constant depending polynomially on $\varepsilon $. Using this fact we obtain
estimates for the covering numbers $N(\sqrt{\smash[b]{q}}B_2^n,tZ_q(\mu ))$, $t\gr 1$, thus showing that $Z_q(\mu )$ is a $\beta $-regular
convex body. As a consequence, we also get an upper bound for $M(Z_q(\mu ))$.
\end{abstract}

\mathversion{bold}
\maketitle
\mathversion{normal}

\section{Introduction}

Given a convex body $K$ of volume 1 or a log-concave probability measure $\mu$
on ${\mathbb R}^n$, we define the $L_q$-centroid bodies $Z_q(K)$ or $Z_q(\mu)$, $q\in (0,+\infty)$,
through their support function $h_{Z_q(K)}$ or $h_{Z_q(\mu)}$, which is defined as follows: for every $y\in {\mathbb R}^n$,
\begin{multline} h_{Z_q(K)}(y):= \|\langle \cdot ,y\rangle \|_{L_q(K)} = \left(\int_K|\langle x,y\rangle|^qdx\right)^{1/q},
\\
h_{Z_q(\mu)}(y):= \|\langle \cdot ,y\rangle \|_{L_q(\mu)} = \left(\int_{{\mathbb R}^n}|\langle x,y\rangle|^qd\mu(x)\right)^{1/q}.
\end{multline}
These bodies then incorporate information about the distribution
of linear functionals with respect to the uniform measure on $K$ or with respect to the probability
measure $\mu$. The $L_q$-centroid bodies were introduced, under a different normalization,
by Lutwak, Yang and Zhang in \cite{LYZ}, while in \cite{PaourisGAFA} for the first time,
and in \cite{PaourisTAMS} later on, Paouris used geometric properties
of them to acquire detailed information about the distribution of the Euclidean norm with
respect to the uniform measure on isotropic convex bodies.
An asymptotic theory for the $L_q$-centroid bodies has since been developed
in the context of isotropic measures and it seems to advance in parallel with all recent developments in the area.

Recall that a convex body $K$ in ${\mathbb
R}^n$ is called isotropic if it has volume $1$, it is centered, i.e.~its
barycenter is at the origin,
and if its inertia matrix is a multiple of the identity matrix:
there exists a constant $L_K >0$ such that
\begin{equation}\label{isotropic-condition}\int_K\langle x,\theta\rangle^2dx =L_K^2\end{equation}
for every $\theta $ in the Euclidean unit sphere $S^{n-1}$.
Similarly, a log-concave probability measure $\mu$ on ${\mathbb R}^n$ is called isotropic
if its barycenter is at the origin and if its inertia matrix is the identity matrix; in that case the isotropic constant
of the measure is defined as
\begin{equation}L_{\mu} := \sup_{x\in {\mathbb R}^n} \bigl(f_{\mu}(x)\bigr)^{1/n},\end{equation}
where $f_{\mu}$ is the density of $\mu$ with respect to the Lebesgue measure.
One very well-known open question in the theory
of isotropic measures is the hyperplane conjecture, which
asks if there exists an absolute constant $C>0$ such that
\begin{equation}\label{HypCon}L_n:= \max\{ L_K:K\ \hbox{is isotropic in}\ {\mathbb R}^n\}\ls C\end{equation}
for all $n\gr 1$. Bourgain proved in \cite{Bou} that $L_n\ls
c\sqrt[4]{n}\log\! n$, while Klartag \cite{Kl}
obtained the bound $L_n\ls c\sqrt[4]{n}$. A second proof of Klartag's bound
appears in \cite{KM1}.

A motivation for this paper is a recent reduction \cite{GPV-JFA} of the hyperplane conjecture
to the study of geometric properties of the $L_q$-centroid
bodies, and in particular to the study of the parameter
\begin{equation}I_1\bigl(K,Z_q^{\circ}(K)\bigr):= \int_K \|x\|_{Z_q^{\circ}(K)}(x) dx = \int_K h_{Z_q(K)}(x)dx.\end{equation}
The main result of \cite{GPV-JFA} is, in a sense, a continuation of
Bourgain's approach to the problem and, roughly speaking, can be
formulated as follows: Given $q\gr 2$ and $\tfrac{1}{2}\ls s\ls 1$,
an upper bound of the form $I_1(K ,Z_q^{\circ }(K))\ls
C_1q^s\sqrt{n}L_K^2$ for all isotropic convex bodies $K$ in
${\mathbb R}^n$ leads to the estimate
\begin{equation}\label{eq:Ln-main}L_n\ls\frac{C_2\sqrt[4]{n}\log\! n}{q^{\frac{1-s}{2}}}.\end{equation}
Bourgain's estimate may be recovered by choosing $q=2$, however, clarifying the
behaviour of $I_1(K ,Z_q^{\circ }(K ))$ might allow one to use much
larger values of $q$. This behaviour is most naturally related to the geometry of
the bodies $Z_q(K)$, and especially how this geometry is affected by or affects the geometry of the body $K$.
This is not yet fully understood and, in view of (\ref{eq:Ln-main}),
we believe that its deeper study would be very useful.

In Section 3 we give an account of some basic known results for the bodies
$Z_q(K)$ and, more generally, the bodies $Z_q(\mu)$
where $\mu$ is an isotropic log-concave measure on ${\mathbb R}^n$.
In the range $2\ls q\ls\sqrt{n}$, for example, some of their global parameters
are completely determined: the volume radius and the mean width of $Z_q(\mu)$ are of the same order $\sqrt{q}$.
The purpose of this work is to provide new information on the local structure of $Z_q(\mu)$,
which in turn has some interesting consequences. Our first main result
concerns proportional projections of the centroid bodies.

\begin{theorem}\label{theorem1.1}Let $\mu $ be an isotropic log-concave measure
on ${\mathbb R}^n$. Fix $1\ls\alpha <2$. For every $0<\varepsilon <1$
and any $q\ls\sqrt{\varepsilon n}$ there are $k\gr (1-\varepsilon
)n$ and $F\in G_{n,k}$ such that
\begin{equation}\label{eq:theorem1.1-eq1}P_F\bigl(Z_q(\mu )\bigr)
\supseteq c(2-\alpha )\varepsilon^{\frac{1}{2}+\frac{2}{\alpha
}}\sqrt{q}\,B_F,\end{equation}
where $c>0$ is an absolute constant $($independent of $\alpha$, $\varepsilon$, the
measure $\mu$, $q$ or $n)$. Moreover, for any $2\ls
q\ls \varepsilon n$ there are $k\gr (1-\varepsilon )n$ and $F\in
G_{n,k}$ such that
\begin{equation}\label{eq:theorem1.1-eq2}P_F\bigl(Z_q(\mu )\bigr)\supseteq \frac{c_1(2-\alpha )\varepsilon^{\frac{1}{2}+\frac{2}{\alpha }}}{L_{\varepsilon n}}\sqrt{q}\,B_F
\supseteq \frac{c_2(2-\alpha )\varepsilon^{\frac{1}{4}+\frac{2}{\alpha }}}{\sqrt[4]{n}}\sqrt{q}\,B_F,\end{equation}
where $c_1, c_2>0$ are absolute constants.
\end{theorem}

The proof of Theorem \ref{theorem1.1} is given in Section 5. We use Pisier's theorem on the existence of $\alpha $-regular
ellipsoids for symmetric convex bodies; we combine this with available information on the $L_q$-centroid bodies
as well as results from \cite{Giannopoulos-Milman-1998} concerning the circumradius of proportional sections of $\alpha $-regular
convex bodies, the proofs of which are outlined in Section 4. Let us mention that the dual result
is a direct consequence of the low $M^{\ast }$-estimate, since
the mean width of $Z_q(\mu )$ is known to be of the order of $\sqrt{q}$: if $2\ls q\ls \sqrt{n}$
and if $\varepsilon\in (0,1)$ and $k=(1-\varepsilon )n$, then a
subspace $F\in G_{n,k}$ satisfies
\begin{equation}\label{eq:theorem1.1-eq3}P_F(Z_q^{\circ }\bigl(\mu )\bigr)\supseteq
\frac{c_1\sqrt{\varepsilon }}{\sqrt{q}} B_F\end{equation} with
probability greater than $1-\exp (-c_2\varepsilon n)$, where $c_1,
c_2>0$ are absolute constants.

In Section 6 we discuss bounds for the covering numbers of a Euclidean ball by $Z_q(\mu )$. It was
proved in \cite{Giannopoulos-Paouris-Valettas-2011} and \cite{Giannopoulos-Paouris-Valettas-2012}
that if $\mu$ is an isotropic log-concave measure on
$\mathbb R^n$ then, for any $1\ls q\ls n$ and $t\gr 1$,
\begin{equation}\label{eq:theorem1.2-eq1}\log{N\bigl( Z_q(\mu),c_1t\sqrt{q}B_2^n\bigr)} \ls
c_2\frac{n}{t^2}+c_3\frac{\sqrt{qn}}{t},\end{equation}
where $c_1,c_2,c_3>0$ are absolute constants. Using Theorem \ref{theorem1.1} and an entropy extension result
from \cite{Litvak-Milman-Pajor-Tomczak-2006} we obtain regular entropy estimates for the dual
covering numbers.

\begin{theorem}\label{theorem1.2}Let $\mu $ be an isotropic log-concave measure on $\mathbb R^n$.
Assume $1\ls\alpha <2$.  Then, for any $q\ls\sqrt{n}$ and any
$$1\ls t\ls \min\bigl\{\sqrt{q},\; c_1(2-\alpha)^{-1}(n/q^2)^{\frac{\alpha+4}{2\alpha}}\bigr\}$$
we have
\begin{equation}\label{eq:theorem1.2-eq2}
\log N\bigl(\sqrt{q}B_2^n, tZ_q(\mu)\bigr) \ls c(\alpha)\frac{n}{t^{\frac{2\alpha }{\alpha +4}}} \max\left\{\log \frac{\sqrt{2q}}{t},
\;\log\frac1{(2-\alpha)t}\right\},\end{equation}
where $c(\alpha)\ls C(2-\alpha )^{-2/3}$ and
$c_1, C$ are absolute constants.
Moreover, for any $2\ls q\ls n$ and any
$$1\ls t\ls\min\biggl\{\sqrt{q},\;
c_2 (2-\alpha)^{-1}L_n\left(\frac{n}{q}\right)^{\frac{\alpha+4}{2\alpha}}\biggr\}$$
we have
\begin{equation}\label{eq:theorem1.2-eq3}
\log N\bigl(\sqrt{q}B_2^n, tZ_q(\mu)\bigr) \ls c(\alpha)L_n^{\frac{2\alpha}{\alpha+4}}\, \frac{n}{t^{\frac{2\alpha }{\alpha +4}}}
\,\max\left\{\log\frac{2q}{t^2},\;\log\frac{L_n}{(2-\alpha)t}\right\},\end{equation}
where $c(\alpha)$ is as above
and $c_2$ is an absolute constant.
\end{theorem}

Note that, since $Z_q(\mu )\supseteq B_2^n$, we are interested in bounds for the above covering numbers when
$t$ is in the interval $[1,\sqrt{\smash[b]{q}}]$. An analysis of the restrictions in Theorem \ref{theorem1.2}
shows that, given any $q\ls n^{3/7}$, (\ref{eq:theorem1.2-eq2}) holds true with any $t$ in the ``interesting'' interval,
while the same is true for (\ref{eq:theorem1.2-eq3}) as long as $q\ls \sqrt{L_n}n^{3/4}$.
Although all these estimates are most probably not optimal, we can still conclude that $Z_q(\mu)$,
with $q\ls n^{3/7}$, is a $\beta $-regular convex body in the sense of Pisier's theorem
(for some concrete positive value of $\beta $). As a consequence of this
fact we get an upper bound for the parameter
\begin{equation*}M\bigl(Z_q(\mu )\bigr)=\int_{S^{n-1}}\| x\|_{Z_q(\mu )}\,d\sigma (x).\end{equation*}
Recall that the dual Sudakov inequality of
Pajor and Tomczak-Jaegermann (see e.g.~\cite{Pisier-book})
provides $2$-regular entropy estimates for the numbers $N(B_2^n,tC)$ in terms of $M(C)$, namely it shows that
\begin{equation*}\log N(B_2^n,tC)\ls cn\left (\frac{M(C)}{t}\right )^2\end{equation*}
for every $t\gr 1$. In Section 7 we use in a converse manner the entropy estimates of Theorem 1.2 to obtain
non-trivial upper bounds for $M(Z_q(\mu ))$.

\begin{theorem}\label{theorem1.3} Let $\mu $ be an isotropic log-concave measure on $\mathbb R^n$.
For every $1\ls q\ls n^{3/7}$,
\begin{equation}\label{eq:theorem1.3-eq1} M\bigl(Z_q(\mu )\bigr)\ls C\frac{(\log q)^{5/6}}{\sqrt[6]{q}}.\end{equation}
Moreover, for every $q$ such that $L_n^2\log^2 q\ls q\ls \sqrt{L_n}\,n^{3/4}$,
\begin{equation}\label{eq:theorem1.3-eq2} M\bigl(Z_q(\mu )\bigr)\ls C\frac{\sqrt[3]{L_n}(\log q)^{5/6}}{\sqrt[6]{q}}.\end{equation}
\end{theorem}

Observe now that, if $K$ is an isotropic convex body in ${\mathbb R}^n$ with isotropic
constant $L_K$,
then the measure $\mu_K$ with density $f_{\mu_K}(x):= L_K^n \mathbf{1}_{{K/L_K}}(x)$ is isotropic
and, for every $q>0$, it holds that $Z_q(K) = L_KZ_q(\mu_K)$. Using also the fact that $M(K)\ls M(Z_q(K))$ for
every symmetric convex body $K$ and every $q>0$, we can use
the above bounds for $M(Z_q(\mu_K))$ to obtain an upper bound for $M(K)$ in the isotropic case.

\begin{theorem}\label{theorem1.4}Let $K$ be a symmetric isotropic convex body in ${\mathbb R}^n$. Then,
\begin{equation*}M(K) \ls C\frac{\sqrt[4]{L_n} (\log n)^{5/6}}{L_K\sqrt[8]{n}}.\end{equation*}
\end{theorem}

This is a question that until recently had not attracted much attention.
Valettas, using a slightly different approach \cite{Valettas-private}, has shown that
\begin{equation*}M(K)\ls \frac{C(\log n)^{1/3}}{\sqrt[12]{n}L_K}\end{equation*}
for every symmetric isotropic convex
body $K$ in ${\mathbb R}^n$, where $C>0$ is an absolute constant.
Note that, on the other hand, there are many approaches concerning the corresponding question about the mean width
that give the best currently known estimate:
\begin{equation*}w(K)\ls Cn^{3/4}L_K\end{equation*}
for every isotropic convex body $K$ in ${\mathbb R}^n$. Nevertheless, this problem as well remains open
(for a discussion about it, see \cite{Giannopoulos-Paouris-Valettas-2012} and the references therein).

We close this paper with some additional observations on the
geometry of the centroid bodies $Z_q(\mu )$ and their polars. We
first provide lower bounds for the radius of their sections;
actually, they hold true for every $1\ls k<n$ and any $F\in
G_{n,k}$. By duality, these estimates (combined with
e.g.~Proposition~\ref{prop:projZq-small}) determine the inradius of
their random projections. We also provide upper bounds for the
parameters $M_{-k}(Z_q(\mu ))$ and $I_{-k}(\overline{Z_q}(\mu ))$
(see Section 8 for the precise definitions). These imply small ball
probability estimates for the Euclidean norm on $Z_q(\mu )$. All the
results are based on the main estimates from Section 5 and Section
6. Although they are not optimal, we describe our approach and
sketch their proofs; we expect that further progress can be achieved
along the same lines.

\section{Notation and preliminaries}

We work in ${\mathbb R}^n$, which is equipped with a Euclidean
structure $\langle\cdot ,\cdot\rangle $. We denote the corresponding
Euclidean norm by $\|\cdot \|_2$, and write $B_2^n$ for the
Euclidean unit ball, and $S^{n-1}$ for the unit sphere. Volume is
denoted by $|\cdot |$. We write $\omega_n$ for the volume of $B_2^n$
and $\sigma $ for the rotationally invariant probability measure on
$S^{n-1}$. We also denote the Haar measure on $O(n)$ by $\nu $. The
Grassmann manifold $G_{n,k}$ of $k$-dimensional subspaces of
${\mathbb R}^n$ is equipped with the Haar probability measure
$\nu_{n,k}$. Let $k\ls n$ and $F\in G_{n,k}$. We will denote the
orthogonal projection from $\mathbb R^{n}$ onto $F$ by $P_F$.
We also define $B_F:=B_2^n\cap F$ and $S_F:=S^{n-1}\cap
F$.

The letters $c,c^{\prime }, c_1, c_2$ etc. denote absolute positive
constants whose value may change from line to line. Whenever we
write $a\simeq b$, we mean that there exist absolute constants
$c_1,c_2>0$ such that $c_1a\ls b\ls c_2a$.  Also if $K,L\subseteq
\mathbb R^n$ we will write $K\simeq L$ if there exist absolute
constants $c_1, c_2>0$ such that $c_{1}K\subseteq L \subseteq
c_{2}K$.

\medskip

\noindent \textbf{Convex bodies.}
A convex body in ${\mathbb R}^n$ is a compact convex subset $C$ of
${\mathbb R}^n$ with nonempty interior. We say that $C$ is symmetric
if $x\in C$ implies that $-x\in C$. We say that $C$ is centered if
the barycenter of $C$ is at the origin, i.e.~$\int_C\langle
x,\theta\rangle \,d x=0$ for every $\theta\in S^{n-1}$. The volume radius of $C$
is the quantity
\begin{equation}\textrm{vrad}(C)=\left (\frac{|C|}{|B_2^n|}\right )^{1/n}.\end{equation}The support
function of a convex body $C$ is defined by
\begin{equation}h_C(y):=\max \bigl\{\langle x,y\rangle :x\in C\bigr\}\end{equation}
and characterizes $C$ uniquely. The mean width of $C$ is
\begin{equation}w(C):=\int_{S^{n-1}}h_C(\theta )\sigma (d\theta )\end{equation}
and the radius of $C$ is the quantity $R(C):=\max\{ \| x\|_2:x\in C\}$. Also, if the origin is an interior point of $C$, the polar body
$C^{\circ }$ of $C$ is defined as follows:
\begin{equation}
C^{\circ}:=\bigl\{y\in {\mathbb R}^n: \langle x,y\rangle \ls 1\;\hbox{for all}\; x\in C\bigr\}.
\end{equation}

\smallskip

Let $C\subset {\mathbb R}^n$ be a symmetric convex body. We write $\|\cdot\|_C$ for the norm
$\| x\|_C=\min\{ t\gr 0:x\in tC\}$ induced to ${\mathbb R}^n$ by $C$, and
we define
\begin{equation}M(C):=\int_{S^{n-1}}\|\theta\|_C d\sigma (\theta ).\end{equation}
Note that $M(C)=w(C^{\circ })$ and that
\begin{equation}M(C)^{-1}\ls \textrm{vrad}(C)\ls w(C)=M(C^{\circ });\end{equation}
the left hand side inequality is easily checked if we express the volume of $C$ as an integral
in polar coordinates and use H\"{o}lder's and Jensen's inequalities, while the right hand side inequality is
the classical Urysohn's inequality. We also need Milman's low $M^{\ast }$-estimate (see \cite{Milman-Schechtman-book}
or \cite{Pisier-book} for precise references):
if $C$ is a symmetric convex body in ${\mathbb R}^n$, then a subspace
$F\in G_{n,k}$ satisfies
\begin{equation}\label{eq:lowM*}R(C\cap F)\ls
c_1\sqrt{\frac{n}{n-k}}\,w(C)\end{equation}
with probability greater than $1-\exp (-c_2(n-k))$, where $c_1, c_2>0$ are absolute
constants.

Furthermore, if $C$ is a symmetric convex body in ${\mathbb R}^n$, we define $k_{\ast }(C)$ to be
the largest positive integer $k\ls n$ with the property that the measure $\nu_{n,k}$ of $F\in G_{n,k}$ for which we
have $\frac{1}{2} w(C)B_F \subseteq P_F(C) \subseteq 2 w(C)B_F$ is
greater than $\frac{n}{n+k}$. It is known that this parameter is completely determined by the dimension,
the mean width and the radius of $C$:
\begin{equation}c_3n\frac{w(C)^2}{R(C)^2} \ls  k_{\ast }(C)\ls  c_4n \frac{w(C)^2}{R(C)^2}.\end{equation}


Recall that the covering number $N(A,B)$ of a body $A$ by a second
body $B$ is the least integer $N$ for which there exist $N$
translates of $B$ whose union covers $A$. Milman
(see e.g.~\cite{Milman-1988}) proved that there exists an absolute constant
$\beta>0$ such that every centered convex body $K$ in ${\mathbb R}^n$ has
a linear image $\tilde{K}$ which satisfies $|\tilde{K}|=|B_2^n|$ and
\begin{equation}\label{betaMposition}
\max\bigl\{ N(\tilde{K},B_2^n), N(B_2^n,\tilde{K}),N(\tilde{K}^{\circ
},B_2^n), N(B_2^n,\tilde{K}^{\circ })\bigr\} \ls\exp(\beta n).
\end{equation}
We say that a convex body $K$ which satisfies this estimate
is in $M$-position with constant $\beta $.

Pisier \cite{Pisier-1989} has proposed a different approach to this result, which
allows one to find a whole family of $M$-ellipsoids and to give
more detailed information on the behaviour of the corresponding covering numbers.
The precise statement is as follows.

\begin{theorem}[Pisier]\label{th:pisier-alpha-regular}For every $0<\alpha <2$
and every symmetric convex body $K$ in ${\mathbb R}^n$ there exists
a linear image $\tilde{K}$ of $K$ such that
\begin{equation*}\max\bigl\{ N(\tilde{K},tB_2^n),N(B_2^n,t\tilde{K}),
N(\tilde{K}^{\circ },tB_2^n), N(B_2^n,t\tilde{K}^{\circ })\bigr\}\ls\exp
\left (\frac{c(\alpha )n}{t^{\alpha }}\right )\end{equation*} for
every $t\gr 1$, where $c(\alpha )$ depends only on $\alpha $, and
$c(\alpha )=O\big ((2-\alpha )^{-\alpha/2}\big )$ as $\alpha\to 2$.
\end{theorem}

For basic facts from the Brunn-Minkowski theory and the asymptotic
theory of finite dimensional normed spaces, whose unit balls are
various symmetric convex bodies appearing in this paper, we refer to
the books \cite{Schneider-book}, \cite{Milman-Schechtman-book} and
\cite{Pisier-book}.

\smallskip

\noindent \textbf{Log-concave probability measures.}
We denote by ${\mathcal{P}}_n$ the class of all Borel
probability measures on $\mathbb R^n$ which are absolutely
continuous with respect to the Lebesgue measure. The density of $\mu
\in {\mathcal{P}}_n$ is denoted by $f_{\mu}$. We say that $\mu
\in {\mathcal{P}}_n$ is centered and we write $\textrm{bar}(\mu )=0$
if, for all $\theta\in S^{n-1}$,
\begin{equation}
\int_{\mathbb R^n} \langle x, \theta \rangle d\mu(x) = \int_{\mathbb
R^n} \langle x, \theta \rangle f_{\mu}(x) dx = 0.
\end{equation}A measure $\mu$ on
$\mathbb R^n$ is called $\log$-concave if $\mu(\lambda
A+(1-\lambda)B) \gr \mu(A)^{\lambda}\mu(B)^{1-\lambda}$ for any compact subsets $A$
and $B$ of ${\mathbb R}^n$ and any $\lambda \in (0,1)$. A function
$f:\mathbb R^n \rightarrow [0,\infty)$ is called $\log$-concave if
its support $\{f>0\}$ is a convex set and the restriction of $\log{f}$ to it is concave.
It is known that if a probability measure $\mu $ is log-concave and $\mu (H)<1$ for every
hyperplane $H$, then $\mu \in {\mathcal{P}}_n$ and its density
$f_{\mu}$ is $\log$-concave. Note that if $K$ is a convex body in
$\mathbb R^n$ then the Brunn-Minkowski inequality implies that
$\mathbf{1}_{K} $ is the density of a $\log$-concave measure.

If $\mu $ is a $\log $-concave measure on ${\mathbb R}^n$ with density $f_{\mu}$,
we define the isotropic constant of $\mu $ by
\begin{equation}\label{definition-isotropic}
L_{\mu }:=\left (\frac{\sup_{x\in {\mathbb R}^n} f_{\mu} (x)}{\int_{{\mathbb
R}^n}f_{\mu}(x)dx}\right )^{\frac{1}{n}} [\det \textrm{Cov}(\mu)]^{\frac{1}{2n}},\end{equation} where
$\textrm{Cov}(\mu)$ is
the covariance matrix of $\mu$ with entries
\begin{equation}\textrm{Cov}(\mu )_{ij}:=\frac{\int_{{\mathbb R}^n}x_ix_j f_{\mu}
(x)\,dx}{\int_{{\mathbb R}^n} f_{\mu} (x)\,dx}-\frac{\int_{{\mathbb
R}^n}x_i f_{\mu} (x)\,dx}{\int_{{\mathbb R}^n} f_{\mu}
(x)\,dx}\frac{\int_{{\mathbb R}^n}x_j f_{\mu}
(x)\,dx}{\int_{{\mathbb R}^n} f_{\mu} (x)\,dx}.\end{equation} We say
that a $\log $-concave probability measure $\mu $ on ${\mathbb R}^n$
is isotropic if $\textrm{bar}(\mu )=0$ and $\textrm{Cov}(\mu )$ is the identity matrix and
we write $\mathcal{IL}_n$ for the class of isotropic $\log $-concave probability measures on ${\mathbb R}^n$.
Note that a centered convex body $K$ of volume $1$ in ${\mathbb R}^n$ is isotropic,
i.e.~it satisfies (\ref{isotropic-condition}),
if and only if the log-concave probability measure $\mu_K$ with density
$x\mapsto L_K^n\mathbf{1}_{K/L_K}(x)$ is isotropic.

Let $\mu\in {\mathcal P}_n$. For every $1\ls k\ls n-1$ and every
$E\in G_{n,k}$, the marginal of $\mu$ with respect to $E$ is the probability
measure with density
\begin{equation}\label{definitionmarginal}f_{\pi_E\mu }(x)= \int_{x+
E^{\perp}} f_{\mu }(y) dy.
\end{equation}
It is easily checked that if $\mu $ is centered, isotropic or log-concave, then $\pi_E\mu $ is also centered, isotropic or
log-concave, respectively.

For more information on isotropic convex bodies and log-concave measures see \cite{Milman-Pajor-1989}, \cite{Ball} and \cite{Gian}.

\section{$L_q$-centroid bodies: basic facts}

Recall that, if $\mu$ is a log-concave probability measure on ${\mathbb R}^n$,
the $L_q$-centroid body $Z_q(\mu )$, $q\gr 1$, of $\mu$ is the centrally symmetric convex body with support
function
\begin{equation}h_{Z_q(\mu )}(y):= \left(\int_{{\mathbb R}^n} |\langle x,y\rangle|^{q}d\mu (x) \right)^{1/q}.\end{equation}
Observe that $\mu $ is isotropic if and only if it is centered and
$Z_2(\mu )=B_2^n$. From H\"{o}lder's inequality it follows that
$Z_1(\mu )\subseteq Z_p(\mu )\subseteq Z_q(\mu )$ for all $1\ls p\ls
q<\infty $. Conversely, using Borell's lemma (see \cite[Appendix III]{Milman-Schechtman-book}), one
can check that
\begin{equation}\label{reverse inclusion for Zq} Z_q(\mu )\subseteq c\frac{q}{p}Z_p(\mu )\end{equation}
for all $1\ls p<q$. In particular, if $\mu $ is isotropic, then $R(Z_q(\mu ))\ls cq$.

As we saw in the previous section, if $K$ is a convex body of volume $1$ in ${\mathbb R}^n$, then the measure
with density $x\mapsto \mathbf{1}_K(x)$ is a log-concave probability measure on ${\mathbb R}^n$
and the $L_q$-centroid bodies $Z_q(K)$ of $K$ can be defined as above. By H\"{o}lder's inequality
we again have $Z_1(K)\subseteq Z_p(K)\subseteq Z_q(K)\subseteq \textrm{conv}\{K,-K\}$ for all $1\ls p\ls
q<\infty $, so, if $K$ is symmetric, then $Z_p(K)\subseteq K$ for every $p<\infty$.

Using Fubini's theorem we see that, for every $1\ls k\ls n-1$ and every $F\in G_{n,k}$ and $q\gr 1$,
\begin{equation}\label{marginalZp}
P_F\bigl(Z_q(\mu )\bigr) = Z_q\bigl(\pi_F(\mu )\bigr).
\end{equation}
In \cite{PaourisGAFA} Paouris shows that the moments
\begin{equation}
I_q(\mu):= \left(\int_{{\mathbb R}^n} \|x\|_2^qdx\right)^{1/q}, \quad q\in (-n,+\infty)\setminus\{0\},
\end{equation}
of the Euclidean norm with respect to an isotropic log-concave probability measure $\mu$ on ${\mathbb R}^n$
remain comparable to $I_2(\mu) = \sqrt{n}$ when $q\gr 2$ does not exceed a parameter $q_{\ast}(\mu)$
of the measure that he defines as follows:
\begin{equation}
q_{\ast}(\mu) := \max\bigl\{ q\ls n : k_{\ast}\bigl(Z_q(\mu )\bigr) \gr q\bigr\}.
\end{equation}
He proves this by showing that
\begin{equation*}I_q(\mu) \simeq \sqrt{\frac{n}{q}} w\bigl(Z_q(\mu)\bigr)\end{equation*}
and that
\begin{equation}\label{eq:wZq-small} w\bigl(Z_q(\mu)\bigr)\simeq \sqrt{q} \end{equation}
for every $q\ls q_{\ast}(\mu)$. He establishes that
\begin{equation}\label{bound-q-ast} q_{\ast}(\mu)\gr c_1\sqrt{n}\end{equation} for every isotropic measure $\mu$ on
${\mathbb R}^n$, where $c_1>0$ is an absolute constant,
and it is also shown that there are isotropic measures $\mu$ on ${\mathbb R}^n$ such that
$q_{\ast}(\mu)\simeq \sqrt{n}$.
Furthemore, he proves that
\begin{equation}I_q(\mu) \simeq R\bigl(Z_q(\mu)\bigr)\end{equation}
for every $q\in [q_{\ast}(\mu), n]$ and he gives an upper bound for the volume radius of the $L_q$-centroid bodies:
\begin{equation}\label{vradZq-upper} |Z_q(\mu)|^{1/n} \ls c_2\sqrt{q/n}\end{equation}
for all $1\ls q\ls n$, where $c_2$ is an absolute constant.

In \cite{PaourisTAMS} Paouris extends his approach to describe the behaviour of the negative moments
of the Euclidean norm with respect to an isotropic measure $\mu$ on ${\mathbb R}^n$. He shows that
\begin{equation}\label{negative-moments} I_{-q}(\mu) \simeq I_2(\mu) = \sqrt{n} \quad\hbox{for all}\  0< q \ls q_{\ast}(\mu).\end{equation}
However, unlike the positive moments $I_{q}(\mu)$ that, as we saw,
do not remain comparable to $I_2(\mu)$ once $q$ gets larger than $q_{\ast}(\mu)$, the behaviour of the corresponding
negative moments is not known and, in fact, (\ref{negative-moments}) may hold with any positive $q$ up to $n-1$.
This question actually is equivalent to the hyperplane conjecture, as Dafnis and Paouris proved in \cite{DP}
by introducing another parameter, that for each $\delta \gr 1$ is given by
\begin{equation}q_{-c}(\mu,\delta):= \max\bigl\{1\ls q\ls n-1 : I_{-q}(\mu)\gr \delta^{-1} I_2(\mu)=
\delta^{-1} \sqrt{n} \bigr\},\end{equation}
namely measures how large the range of $(\ref{negative-moments})$ is
if we allow the implied constants to depend on $\delta$.
Dafnis and Paouris established that
\begin{equation}\label{eq:mainDP}
L_n\ls C\delta\sup_{\mu\in {\mathcal IL}_{[n]}}\sqrt{\frac{n}{q_{-c}(\mu,\delta)}}\;
\log\Bigl(\frac{en}{q_{-c}(\mu,\delta)}\Bigr)
\end{equation}
for every $\delta\gr 1$, and they also showed that, if the hyperplane conjecture is correct,
that is, if (\ref{HypCon}) holds true, then we will have
\begin{equation} q_{-c}\bigl(\mu,\delta_0\bigr)=n-1 \end{equation}
for some $\delta_0 \simeq 1$, for every isotropic log-concave measure $\mu$ on ${\mathbb R}^n$.
Note that, by (\ref{negative-moments}), we already know that
\begin{equation}\label{bound-q--c}  q_{-c}\bigl(\mu,\delta_1\bigr) \gr q_{\ast}(\mu) \gr c_1\sqrt{n},\end{equation}
where $\delta_1\gr 1$ and $c_1>0$ are absolute constants.

Next we turn to lower bounds for the volume radius of
the $L_q$-centroid bodies of measures $\mu\in \mathcal{IL}_{[n]}$.
From \cite{LYZ} we know that, for every $1\ls q\ls n$,
\begin{equation}\label{eq0:low-volume-Zq}
\bigl|Z_q(\mu)\bigr|^{1/n} \gr c_2\sqrt{q/n}\, L_{\mu}^{-1}
\end{equation}
for some absolute constant $c_2>0$. In \cite{KM1} Klartag and Milman define a ``hereditary" variant of $q_{\ast }(\mu )$ as follows:
\begin{equation}q_{\ast}^H (\mu) := n\inf_k \inf_{E\in G_{n,k}}\frac{q_{\ast}(\pi_E\mu)}{k},\end{equation}
where $\pi_E\mu $ is the marginal of $\mu $ with respect to $E$, and then, for every $q\ls q_{\ast}^H (\mu)$, they
give a lower bound for the volume radius of the bodies $Z_q(\mu)$ matching the upper bound in (\ref{vradZq-upper}):
\begin{equation}\label{eq:low-volume-Zq}
\bigl|Z_q(\mu)\bigr|^{1/n}\gr c_3\sqrt{q/n}
\end{equation}
where $c_3>0$ is an absolute constant. Recall that, if $\mu$ is an isotropic log-concave measure, then so are all its marginals,
thus, for every subspace $E\in G_{n,k}$, we have by (\ref{bound-q-ast}) that $q_{\ast}(\pi_E\mu)\gr c_1\sqrt{k}$.
This implies that, for all measures $\mu\in \mathcal{IL}_{[n]}$, $q_{\ast}^H (\mu)\gr c_1\sqrt{n}$. Note also that,
for those measures for which we have $q_{\ast} (\mu)\simeq \sqrt{n}$, the parameter $q_{\ast}^H (\mu)$ as well
does not exceed a constant multiple of $\sqrt{n}$. However, the bound (\ref{eq:low-volume-Zq}) might hold for larger $q\in [1,n]$
even in the latter case, as was shown by the fourth named author in \cite{Vritsiou-IMRN}. There a hereditary variant
of $q_{-c}(\mu, \delta)$ is introduced as follows:
\begin{equation*}
q_{-c}^H(\mu,\delta):= n \inf_k\inf_{E\in G_{n,k}} \frac{ q_{-c}(\pi_E\mu,\delta)}{k}
\end{equation*}
for every $\delta \gr 1$, and in a similar way as for the bound (\ref{eq:low-volume-Zq}) it is established that
\begin{equation}\label{eq1:low-volume-Zq}
\bigl|Z_q(\mu)\bigr|^{1/n} \gr c_4 \delta^{-1}\sqrt{q/n}
\end{equation}
for every $q\ls q_{-c}^H(\mu,\delta)$. Of course, we can again see that
\begin{equation}q_{-c}^H(\mu,\delta_1)\gr c_1\sqrt{n}\end{equation}
for some $\delta_1\simeq 1$ using (\ref{bound-q--c}) and the definition of $q_{-c}^H(\mu,\delta_1)$,
but obviously this estimate might be improved if the latter bound
were also; such an improvement may be possible but is not trivial, since it would result in better bounds for the isotropic constant
problem too, as one can see from (\ref{eq:mainDP}).

To conclude this Section, we should stress that in the subsequent
proofs we will be using the lower bounds for the volume radius of
the $L_q$-centroid bodies, however the only concrete estimate which
we currently have for the parameters $q_{\ast }^H(\mu )$ or
$q_{-c}^H(\mu ,\delta )$, and which we can insert in our later
estimates, is that all of them are at least of the order of
$\sqrt{n}$ when $\mu$ is an isotropic log-concave measure on
${\mathbb R}^n$. For more information on all the results of this
Section, see \cite{BGVV-book-isotropic}.

\section{Diameter of sections of $\alpha $-regular bodies}

We say that a body $\tilde{K}$ which satisfies the conclusion of
Theorem \ref{th:pisier-alpha-regular} is an \textit{$\alpha $-regular
body}. A strong form of the reverse Brunn-Minkowski inequality can
be proved for $\alpha $-regular bodies; actually, we only need the
regularity of the covering numbers $N(K,tB_2^n)$.

\begin{lemma}\label{th:M-7}
Let $\gamma\gr 1$, $\alpha >0$ and let $K_1,\ldots ,K_m$ be
symmetric convex bodies in ${\mathbb R}^n$ which satisfy
\begin{equation*}N(K_j,tB_2^n)\ls \exp\left (\frac{\gamma
n}{t^{\alpha }}\right)\end{equation*}for all $1\ls j\ls m$ and all $t\gr 1$. Then,
\begin{equation}|K_1+\cdots +K_m|^{1/n}\ls
C\gamma^{\frac{1}{\alpha }}m^{1+\frac{1}{\alpha}
}|B_2^n|^{1/n}.\end{equation}
\end{lemma}

\noindent \textit{Proof.} We include the very simple proof of this fact. Observe that
\begin{align*}N(K_1+\cdots +K_m,tmB_2^n)&=N(K_1+\cdots
+K_m,tB_2^n+\cdots +tB_2^n) \\
&\ls\prod_{j=1}^mN(K_j,tB_2^n)\ls\exp (\gamma nm/t^{\alpha
})\end{align*}for all $t\gr 1$. It follows that
\begin{equation*}
|K_1+\cdots +K_m|^{1/n} \ls tm\,\exp (\gamma m/t^{\alpha
})|B_2^n|^{1/n}.
\end{equation*}
Choosing $t=(\gamma m)^{1/\alpha }$ we get the result. $\hfill\Box $

\bigskip

We also need to recall the precise probabilistic form of the low $M^{\ast
}$-estimate (which can be found in \cite{Gordon} and \cite{Milman-1991}): If $A$ is a
symmetric convex body in ${\mathbb R}^n$ and if $\varepsilon
,\delta\in (0,1)$, then we have \begin{equation*}R(A\cap
F)\leq\frac{w(A)}{(1-\delta ) \sqrt{\varepsilon}}\end{equation*} for
all $F$ in a subset $L_{n,k}$ of $G_{n,k}$ of measure
$\nu_{n,k}(L_{n,k})\gr 1-c_1\exp (-c_2\delta^2\varepsilon n)$, where
$k=\lfloor (1-\varepsilon )n\rfloor $ and $c_1,c_2>0$ are absolute
constants. Then, a well-known application of the low $M^{\ast
}$-estimate (see e.g.~\cite[Theorem 2.1]{Giannopoulos-Milman-1997})
states that if $r>0$ is the solution of the equation
\begin{equation}\label{eq:M*-equation}\frac{w(A\cap
rB_2^n)}{r}=\frac{1}{2}\sqrt{\varepsilon },
\end{equation} then a typical $\lfloor (1-\varepsilon )n\rfloor$-dimensional
central section of $A$ has radius smaller than $r$ (with probability
greater than $1-\exp (-c\varepsilon n)$).

\medskip

The main observation of this Section, which essentially appears in \cite{Giannopoulos-Milman-1998}, is the following:
if $A$ is a symmetric convex body in ${\mathbb R}^n$ whose covering numbers $N(A,tB_2^n)$
are $\alpha $-regular for some $\alpha >0$, then one can get an upper bound for the diameter of random
proportional sections of $A$.

\begin{theorem}\label{th:section-alpha-regular}Let $\gamma\gr 1$, $\alpha >0$ and let $A$
be a symmetric convex body in ${\mathbb R}^n$ which satisfies
\begin{equation*}N(A,tB_2^n)\ls \exp\left (\frac{\gamma n}{t^{\alpha
}}\right)\end{equation*}for all $t\gr 1$. Then, for every
$\varepsilon\in (0,1)$ a subspace $F\in G_{n,\lfloor
(1-\varepsilon ) n\rfloor }$ satisfies
\begin{equation*}R(A\cap F)\ls C\gamma^{\frac{1}{\alpha }}/\varepsilon^{\frac{1}{2}+\frac{1}{\alpha }},\end{equation*}
with probability greater than $1-c_1\exp (-c_2\varepsilon n)$, where $c_1,c_2,C>0$ are
absolute constants.\end{theorem}

\noindent \textit{Proof.} Let $\varepsilon\in (0,1)$ and set $k=\lfloor
(1-\varepsilon )n\rfloor $. We define $r>0$ by the equation
\begin{equation*}w(A\cap rB_2^n)=\frac{1}{2}\sqrt{\varepsilon }\,r.\end{equation*}
From the precise probabilistic form of the low $M^{\ast}$-estimate (see above), we know that there exists a subset
$L_{n,k}$ of $G_{n,k}$ with measure $\nu_{n,k}(L_{n,k})\gr 1-c_1\exp
(-c_2\varepsilon n)$, such that
\begin{equation*}w(A\cap F)\ls R(A\cap F
)\ls r\end{equation*} for every $F\in L_{n,k}$. We use the following
fact from \cite{BLM}: if $X=({\mathbb R}^n,\|\cdot\|)$ is an
$n$-dimensional normed space with unit ball $W$, and if
$M=\int_{S^{n-1}}\| x\|\,d\sigma (x)$ and $b$ is the smallest
positive constant for which $\| x\|\ls b\|x\|_2$ for all $x\in
{\mathbb R}^n$, then there exist an integer $s\ls C(b/M)^2$ and $s$
orthogonal transformations $U_1,\ldots ,U_s\in O(n)$ such that
\begin{equation*}\frac{M}{2}\,B_2^n\subseteq
\frac{1}{s}\sum_{i=1}^sU_i(W^{\circ })\subseteq 2M\,B_2^n.\end{equation*}
We apply this result for the body $W:=(A\cap rB_2^n)^{\circ }$. Note that $b=r$ and $M(W)=w(A\cap rB_2^n)=
\sqrt{\varepsilon }r/2$, and hence we can find $s\ls\frac{c_3}{\varepsilon }$ and
orthogonal transformations $U_1,\ldots ,U_s$, satisfying
\begin{equation}\label{eq:section-alpha-regular-1}\frac{1}{4}\sqrt{\varepsilon }\,rB_2^n\subseteq\frac{1}{s}
\sum_{i=1}^sU_i(A\cap rB_2^n)\subseteq \sqrt{\varepsilon
}\,rB_2^n.\end{equation} Set $A_1=\frac{1}{s}\sum U_i(A\cap
rB_2^n)$. Now we can give an upper bound for $r$ using
Lemma~\ref{th:M-7}. Clearly, the bodies $U_i(A\cap rB_2^n)$ satisfy
\begin{equation*}
N\bigl(U_i(A\cap rB_2^n),tB_2^n\bigr)\ls N(A,tB_2^n)\ls \exp\left
(\frac{\gamma n}{t^{\alpha }}\right)
\end{equation*}
for all $t\gr 1$, therefore
\begin{equation*}\frac{1}{4}\sqrt{\varepsilon }\,r\ls\left
(\frac{|A_1|}{|B_2^n|}\right )^ {\frac{1}{n}}\ls
c_4\gamma^{\frac{1}{\alpha }}s^{\frac{1}{\alpha} }.\end{equation*}
This shows that
\begin{equation*}R(A\cap F)\ls r\ls c_5\gamma^{\frac{1}{\alpha }}/
\varepsilon^{\frac{1}{2}+\frac{1}{\alpha }}\end{equation*}with
probability greater than $1-c_1\exp (-c_2\varepsilon n)$.
$\hfill\Box $

\section{Projections of $L_q$-centroid bodies}

Our aim is to obtain lower bounds for the inradius of proportional
projections of $Z_q(\mu )$ and $Z_q^{\circ }(\mu )$.
Let $1\ls k\ls n-1$ and consider a random subspace $F\in G_{n,k}$.
An upper bound for the radius of $Z_q(\mu )\cap F$, and hence a lower bound for
the inradius of $P_F(Z_q^{\circ }(\mu ))$, follows from the
low $M^{\ast }$-estimate (\ref{eq:lowM*}) and (\ref{eq:wZq-small}).

\begin{proposition}\label{prop:upper-ZqKcapF-small}Let $\mu $ be an isotropic
log-concave measure on ${\mathbb R}^n$. If $2\ls q\ls q_{\ast}(\mu)$
and if $\varepsilon\in (0,1)$ and $k=\lfloor (1-\varepsilon )n\rfloor $, then a
subspace $F\in G_{n,k}$ satisfies
\begin{equation}\label{eq:upper-ZqKcapF-small}
R\bigl(Z_q(\mu )\cap F\bigr)\ls
\frac{c_1\sqrt{q}}{\sqrt{\varepsilon }}\quad\hbox{or equivalently}\quad
P_F\bigl(Z_q^{\circ }(\mu )\bigr)\supseteq \frac{c_2\sqrt{\varepsilon }}{\sqrt{q}}\,B_F\end{equation} with
probability greater than $1-c_3\exp (-c_4\varepsilon n)$, where $c_i$ are absolute constants. $\hfill\Box $
\end{proposition}

We provide analogous upper bounds for $R(Z_q^{\circ
}(\mu )\cap F)$, $F\in G_{n,k}$. The idea of the proof comes from
\cite{KlM1} (see the concluding remarks of this section). We start with the following
immediate consequence of Theorem \ref{th:section-alpha-regular}.

\begin{corollary}\label{cor:projZq-1}Let $A$ be a symmetric
convex body in ${\mathbb R}^m$. Assume that $\gamma\gr 1$, $\alpha >0$, and $E$ is
an ellipsoid in ${\mathbb R}^m$ such that
\begin{equation*}N(A,tE)\ls\exp \left (\frac{\gamma m}{t^{\alpha }}\right )\end{equation*}
for all $t\gr 1$. Then, for every $\varepsilon\in (0,1)$ there
exists $F\in G_{m,\lfloor (1-\varepsilon ) m\rfloor }$ such that
\begin{equation*}A\cap F\subseteq c\gamma^{\frac{1}{\alpha }}
\varepsilon^{-(\frac{1}{2}+\frac{1}{\alpha })}\,E\cap
F,\end{equation*}where $c>0$ is an absolute constant.\end{corollary}

\begin{proposition}[version for ``small'' $q$]\label{prop:projZq-small}
Let $\mu $ be an isotropic log-concave measure on
${\mathbb R}^n$. Let $1\ls\alpha <2$. For every $0<\varepsilon <1$
and any $q\ls\sqrt{\varepsilon n}$ there exist $k\gr (1-\varepsilon
)n$ and $F\in G_{n,k}$ such that
\begin{equation*}P_F\bigl(Z_q(\mu )\bigr)\supseteq c(2-\alpha )\varepsilon^{\frac{1}{2}+\frac{2}{\alpha
}}\sqrt{q}\,B_F.\end{equation*}
\end{proposition}

\noindent \textit{Proof.} Recall from (\ref{marginalZp}) that for every $1\ls m\ls n$ and any
$H\in G_{n,m}$ we have
\begin{equation*}
P_H\bigl(Z_q(\mu )\bigr) = Z_q\bigl(\pi_H(\mu )\bigr).
\end{equation*}
In Section 3 we saw that, if $\nu $ is an isotropic log-concave
measure on ${\mathbb R}^m$, then
\begin{equation*}|Z_q(\nu )|^{1/m}\gr c\sqrt{q/m}\end{equation*}for
all $q\ls q_{-c}^H(\nu ,\delta_0)$, where $\delta_0\gr 1$ is an absolute
constant sufficiently large so that $q_{-c}^H(\nu ,\delta_0)\gr c\sqrt{m}$. It follows that
\begin{equation}\label{eq:projZq-1}\Bigl|P_H\bigl(Z_q(\mu )\bigr)\Bigr|^{1/m}\gr
c_1\sqrt{q/m}\end{equation}for all $H\in G_{n,m}$ and all $q\ls
\sqrt{m}$. We fix $1\ls\alpha <2$ and consider an $\alpha $-regular
$M$-ellipsoid $E$ of $Z_q(\mu )$, namely an ellipsoid such that
\begin{equation*}
\max\Bigl\{N\bigl(Z_q(\mu),tE\bigr),
N\bigl(E,tZ_q(\mu )\bigr)\Bigr\}\ls e^{c(\alpha )n/t^{\alpha }}\end{equation*}
for all $t\gr 1$, where $c(\alpha )\ls C(2-\alpha)^{-\alpha /2}$.

Let $0<\lambda_1\ls\cdots\ls \lambda_n$ be the axes of $E$, and let
$\{u_1,\ldots ,u_n\}$ be an orthonormal basis which corresponds to the
$\lambda_j$. For every $1\ls m,s\ls n$ we
set\begin{equation*}H_m:=\textrm{span}\{ u_1,\ldots
,u_m\}\quad\hbox{and}\quad F_s=\textrm{span}\{ u_{s+1},\ldots
,u_n\}.\end{equation*}Since $E\cap H_m = P_{H_m}(E)$, we
have\begin{equation*}N\Bigl(P_{H_m}\bigl(Z_q(\mu )\bigr),t(E\cap H_m)\Bigr)\ls
N\bigl(Z_q(\mu),tE\bigr)\ls
e^{c(\alpha )n/t^{\alpha }},\end{equation*}and hence
\begin{equation}\label{eq:projZq-2}\Bigl|P_{H_m}\bigl(Z_q(\mu )\bigr)\Bigr|^{1/m}\ls
e^{c(\alpha )n/(t^{\alpha
}m)}|B_{H_m}|^{1/m}(t\lambda_m).\end{equation}Choose $t=(c(\alpha
)n/m)^{1/\alpha }$. Assuming that $q\ls\sqrt{m}$, from
(\ref{eq:projZq-1}) and (\ref{eq:projZq-2}) we
get\begin{equation}\label{eq:lambda}\lambda_m\gr \left
(\frac{m}{c(\alpha )n}\right )^{1/\alpha
}\sqrt{q}.\end{equation}Next, let $0<\varepsilon <1$ and set
$s=\lfloor\frac{\varepsilon n}{2}\rfloor $. We have
\begin{equation*}N\Bigl(E\cap F_s,tP_{F_s}\bigl(Z_q(\mu )\bigr)\Bigr)\ls
N\bigl(E, tZ_q(\mu )\bigr)\ls
e^{c(\alpha )n/t^{\alpha }}\ls e^{2c(\alpha )(n-s)/t^{\alpha}},\end{equation*}for every $t\gr 1$.

We now use the duality of entropy theorem of Artstein-Avidan,
Milman and Szarek~\cite{Artstein-Milman-Szarek-2004}:
There exist two absolute constants $a$ and $b>0$ such that
for any dimension $n$ and any symmetric convex body $A$ in ${\mathbb
R}^n$ one has
\begin{equation*}
{N(B_2^n, a^{-1}A^{\circ})}^{1/b} \ls {N(A, B_2^n)} \ls
{N(B_2^n, aA^{\circ})}^b.
\end{equation*}
It follows that
\begin{equation*}N(Z_q^{\circ }(\mu )\cap F_s, tE^{\circ }\cap
F_s)\ls N\Bigl(E\cap F_s,atP_{F_s}\bigl(Z_q(\mu )\bigr)\Bigr)^b\ls e^{c_1(\alpha
)(n-s)/t^{\alpha }}.\end{equation*} We apply Corollary \ref{cor:projZq-1}
with the body $Z_q^{\circ }(\mu )\cap F_s$ (and $\gamma =c_1(\alpha )$) to
find a subspace $F$ of $F_s$, of dimension $k\gr (1-\varepsilon
/2)(n-s)\gr (1-\varepsilon )n$, such that
\begin{equation*}Z_q^{\circ }(\mu )\cap F\subseteq
\frac{C} {\sqrt{2-\alpha }\,\varepsilon^{\frac{1}{2}+\frac{1}{\alpha
}}}\,E^{\circ }\cap F,\end{equation*}and hence
\begin{equation}\label{eq:proj-Zq-3}P_F\bigl(Z_q(\mu)\bigr)\supseteq c\sqrt{2-\alpha
}\, \varepsilon^{\frac{1}{2}+\frac{1}{\alpha }}
P_F(E).\end{equation}From (\ref{eq:lambda}) we have
\begin{equation*}E\cap F_s\supseteq \lambda_{s+1}B_{F_s}\supseteq
c\sqrt{2-\alpha }\,\varepsilon^{1/\alpha
}\sqrt{q}\,B_{F_s},\end{equation*}provided that
$q\ls\sqrt{\varepsilon n}$. Then,
\begin{align*}P_F(E) &=P_F\bigl(P_{F_s}(E)\bigr)=P_F(E\cap F_s)\supseteq c\sqrt{2-\alpha
}\,\varepsilon^{1/\alpha }\sqrt{q}P_F(B_{F_s})\\
&=c\sqrt{2-\alpha }\,\varepsilon^{1/\alpha }\sqrt{q}B_F.\end{align*}
Combining this fact with \eqref{eq:proj-Zq-3} we conclude the proof.
$\hfill\Box $

\begin{proposition}[version for ``large'' $q$]\label{prop:projZq}
Let $\mu $ be an isotropic log-concave measure on ${\mathbb R}^n$.
Let $1\ls \alpha <2$. For every $0<\varepsilon <1$ and any $2\ls
q\ls \varepsilon n$ there exist $k\gr (1-\varepsilon )n$ and $F\in
G_{n,k}$ such that
\begin{equation*}P_F\bigl(Z_q(\mu )\bigr)\supseteq \frac{c_1(2-\alpha )\varepsilon^{\frac{1}{2}+\frac{2}{\alpha }}}{L_{\varepsilon n}}\sqrt{q}\,B_F
\supseteq \frac{c_2(2-\alpha )\varepsilon^{\frac{1}{4}+\frac{2}{\alpha }}}{\sqrt[4]{n}}\sqrt{q}\,B_F.
\end{equation*}\end{proposition}

\noindent \textit{Proof.} We apply the same argument more or less, only instead of the lower bounds
(\ref{eq:low-volume-Zq}), (\ref{eq1:low-volume-Zq}) for the volume radius of the
$L_q$-centroid bodies we use (\ref{eq0:low-volume-Zq}).
It follows that
\begin{equation}\label{eq:projZq-5}\Bigl|P_H\bigl(Z_q(\mu )\bigr)\Bigr|^{1/m}\gr
\frac{c_1}{L_m}\sqrt{q/m}\end{equation}for all $H\in G_{n,m}$ and
all $q\ls m$. We define $H_m, F_s$ as in the proof of
Proposition~\ref{prop:projZq-small} and we consider an $\alpha
$-regular $M$-ellipsoid $E$ of $Z_q(\mu )$. This time, assuming that $q\ls
m$, from (\ref{eq:projZq-5}) and (\ref{eq:projZq-2}) we
get\begin{equation*}\lambda_m\gr \frac{1}{L_m}\left
(\frac{m}{c(\alpha )n}\right )^{1/\alpha
}\sqrt{q}.\end{equation*}
Next, fix some $\varepsilon \in (0,1)$,
set $s=\lfloor\frac{\varepsilon n}{2}\rfloor $ and consider any
$q\ls\varepsilon n/2$. As previously, we find a
subspace $F$ of $F_s$, of dimension $k\gr (1-\varepsilon /2)(n-s)\gr
(1-\varepsilon )n$, such that
\begin{equation*}P_F\bigl(Z_q(\mu )\bigr) \supseteq c\sqrt{2-\alpha
}\,\varepsilon^{\frac{1}{2}+\frac{1}{\alpha }} P_F(E)\end{equation*}
and
\begin{equation*}P_F(E) = P_F\bigl(P_{F_s}(E)\bigr)\supseteq P_F(\lambda_{s+1}B_{F_s})
\supseteq \frac{c\sqrt{2-\alpha}}{L_s}\varepsilon^{1/\alpha}\sqrt{q}\,B_F.\end{equation*}
Since $s\simeq \varepsilon n$ and $L_{\varepsilon n}\ls
C\sqrt[4]{\varepsilon n}$, the result follows. $\hfill\Box $

\medskip

As we saw in Section 3, if $K$ is an isotropic symmetric convex body in ${\mathbb
R}^n$, then
\begin{equation*}P_F(K)\supseteq P_F\bigl(Z_q(K)\bigr)\end{equation*}for all
$q > 0$. Recall also that the measure $\mu_K$ with density
$L_K^n\mathbf{1}_{{K/L_K}}$ is isotropic and
$Z_q(K)=L_KZ_q(\mu_K)$. Choosing $q=\varepsilon n$ and applying
Proposition \ref{prop:projZq} with $\mu=\mu_K$ we get:

\begin{corollary}\label{cor:projZq(K)}
Let $K$ be a symmetric isotropic convex body in ${\mathbb R}^n$.
For every $1\ls\alpha <2$ and $0<\varepsilon <1$ there exist $k\gr
(1-\varepsilon )n$ and $F\in G_{n,k}$ such that
\begin{equation*}P_F(K)\supseteq L_KP_F\bigl(Z_q(\mu_K)\bigr)
\supseteq c(2-\alpha )\,\varepsilon^{\frac{3}{4}+\frac{2}{\alpha }}\sqrt[4]{n}L_K\,B_F.\end{equation*}
\end{corollary}

\begin{remark}\rm Some variants of Corollary \ref{cor:projZq(K)} have appeared in the
literature before. In \cite{Valettas-private} a stronger estimate is obtained
with a different method: if $K$ is
a symmetric isotropic convex body in $\mathbb R^n$ then, for
any $0<\varepsilon <1$, there exists a subspace $F$ of ${\mathbb
R}^n$ with $\textrm{dim}\,F\gr (1-\varepsilon)n$ such that
\begin{equation*}
P_F(K)\supseteq c \varepsilon^{3/2} \frac{\sqrt[4]{n}}{\log n}L_K
B_F,\end{equation*} where $c>0$ is an absolute constant. A similar result with cubic dependence
on $\varepsilon$ appears in \cite{Bastero-2007}. Our argument is very much related to the one in \cite{KlM1}
where, under the additional assumption that $L_n\ls C$ for all $n\gr 1$, the existence of some $F\in G_{n,\lfloor (1-\varepsilon )n\rfloor }$
so that
\begin{equation*}P_F(K)\supseteq c \varepsilon^3\sqrt{n}B_F\end{equation*}
is established for all isotropic convex bodies $K$ in ${\mathbb R}^n$ and all $0<\varepsilon <1$. Under this assumption,
our argument would result in the estimate
$P_F(K)\supseteq c(2-\alpha )\varepsilon^{\frac{\alpha +2}{\alpha }}\sqrt{n}B_F$ for all $1\ls\alpha <2$.
\end{remark}

\begin{remark}\rm Proposition \ref{prop:projZq-small} and Proposition \ref{prop:projZq} guarantee
the existence of \textit{one} $\lfloor (1-\varepsilon )n\rfloor $-dimensional projection of $Z_q(\mu )$
with ``large" inradius. However, it is proved in \cite{Giannopoulos-Milman-Tsolomitis-2005} that,
for every fixed proportion $\mu\in
(0,1)$ and every $0<s< 1/(2-\mu )$, the maximal inradius of
$\lfloor\mu n\rfloor $-dimensional projections and the random
inradius of $\lfloor s\mu  n\rfloor $-dimensional projections of a symmetric convex body $K$
in ${\mathbb R}^n$ are comparable up to a constant depending on $\mu $ and $s$. More
precisely, if $a(\lambda , K)$ denotes the maximal (and if $b(\lambda
,K)$ denotes the ``random") inradius of a $\lfloor \lambda
n\rfloor $-dimensional projection of $K$ then
\begin{equation*}\left (\frac{c\mu \bigl(1-s(2-\mu )\bigr)}{1-s\mu } \sqrt{1-\mu }\right )a(\mu , K)\ls b(s\mu  , K)
\end{equation*} for every $n\gr n_0(\mu ,s)$. Using this fact one can obtain versions
of the results of this Section concerning random proportional projections of $Z_q(\mu )$. Since an estimate
for the maximal inradius is sufficient for our subsequent work in this paper, we do not present
the precise statements.
\end{remark}

\section{Covering numbers}

Using Proposition \ref{prop:projZq-small}, for any isotropic
log-concave measure $\mu $ on ${\mathbb R}^n$ we can get some
estimates for the covering numbers $N(\sqrt{q}B_2^n,tZ_q(\mu ))$.
These will follow from an entropy extension result from
\cite{Litvak-Milman-Pajor-Tomczak-2006}:

\begin{lemma}\label{lem:covering}
Let $K,L$ be symmetric convex bodies in $\mathbb R^n$ and assume
that $L\subseteq RK$. Let $F$ be a subspace of $\mathbb R^n$ with
$\dim F=n-m$ and let $0<r<t<R$. Then, we have
\begin{equation}\label{eq:covering-1}
N(L,tK)\ls 2^m\left(\frac{2R+t}{t-r}\right)^mN\left(P_F(L),
\frac{r}{2}P_F(K) \right).
\end{equation}
\end{lemma}

\begin{remark}\upshape Alternatively, one might use an analogous result, due to
Vershynin and Rudelson (see \cite[Lemma 5.2]{Vershynin-2006}): If $K$ is a symmetric convex body
in ${\mathbb R}^n$ such that $K\supseteq\delta B_n$ and if $P_F(K)\supseteq B_F$ for some
$F\in G_{n,k}$, $k\gr (1-\varepsilon )n$, then
\begin{equation*}N(B_2^n,4K)\ls (C/\delta )^{2\varepsilon n}.\end{equation*}The reader may wish to
check that applying this fact instead of Lemma~\ref{lem:covering} leads to the
same estimate in Proposition \ref{prop:covering} below.
\end{remark}

\begin{proposition}\label{prop:covering}
Let $\mu $ be an isotropic log-concave measure on $\mathbb R^n$.
Assume that $q\ls\sqrt{n}$. Then, for any $1\ls\alpha <2$ and any
$$1\ls t\ls \min\bigl\{\sqrt{q},\; c_2(2-\alpha)^{-1}(n/q^2)^{\frac{\alpha+4}{2\alpha}}\bigr\}$$ we have
\begin{eqnarray} 
&& \max\Bigl\{ \log N\bigl(\sqrt{q}B_2^n, tZ_q(\mu)\bigr),
\log N\bigl(\sqrt{q}Z_q^{\circ }(\mu),tB_2^n\bigr)\Bigr\}\nonumber\\[1ex]
&&\hspace*{12em}\ls c(\alpha)\frac{n}{t^{\frac{2\alpha }{\alpha +4}}} \max\left\{\log \frac{\sqrt{2q}}{t},
\;\log\frac1{(2-\alpha)t}\right\},\label{eq:covering-2}
\end{eqnarray}
where $c(\alpha )\ls C(2-\alpha )^{-2/3}$.
\end{proposition}

\noindent \textit{Proof.} Note that, since $B_2^n\subseteq Z_q(\mu )$,
the interesting range for $t$ is up to
$\sqrt{q}$. Given some $\varepsilon\in (0,1)$, let
$k=(1-\varepsilon )n$ and $F\in G_{n,k}$. Applying
Lemma~\ref{lem:covering} for the bodies $\sqrt{q}B_2^n$ and $Z_q(\mu )$
with $R=\sqrt{q}$ and $r=t/2$ we see that, for every $1\ls t<\sqrt{q}$,
\begin{equation}\label{eq:covering-3}
N\bigl(\sqrt{q}B_2^n,tZ_q(\mu )\bigr)\ls
\left(\frac{c_1\sqrt{q}}{t}\right)^{\varepsilon n}N\left(
\sqrt{q}B_F,\frac{t}{4}P_F\bigl(Z_q(\mu )\bigr)\right).
\end{equation} If $q^2\ls\varepsilon n$ then Proposition \ref{prop:projZq-small}
shows that, for every $1\ls\alpha <2$, there exists $F\in G_{n,k}$,
$k=(1-\varepsilon )n$, such that $P_F(Z_q(\mu ))\supseteq
c_2(2-\alpha )\varepsilon^{\frac{1}{2}+\frac{2}{\alpha
}}\sqrt{q}B_F$. Thus, we arrive at
\begin{equation}\label{eq:covering-4}N\bigl(\sqrt{q}B_2^n,tZ_q(\mu )\bigr)\ls\left(\frac{c_3\sqrt{q}}{t}\right)^{\varepsilon n}
N\left(B_F, c_4t(2-\alpha )\varepsilon^{\frac{1}{2}+\frac{2}{\alpha
}}B_F\right).
\end{equation} In the end we choose $\varepsilon\simeq [(2-\alpha )t]^{-\frac{2\alpha }{\alpha +4}}$
(the restriction $[(2-\alpha )t]^{\frac{2\alpha }{\alpha
+4}}\ls cn/q^2$ is needed at this point---note that, if e.g.~$q\ls n^{3/7}$,
then this allows us to consider any $t$ up to $\sqrt{q}$). With this choice
of $\varepsilon$, we get from \eqref{eq:covering-4} that
\begin{equation}\label{eq:covering-5}\log N\bigl(\sqrt{q}B_2^n,tZ_q(\mu )\bigr)
\ls c(\alpha )\frac{n\log (2q/t^2)}{t^{\frac{2\alpha }{\alpha +4}}}.
\end{equation}
This proves the upper bound for the first covering number in (\ref{eq:covering-2})
provided that $t\gr c_1(2-\alpha)^{-1}$ (since $(2-\alpha )t\simeq \varepsilon^{\frac{\alpha+4}{2\alpha }}$
must be less than 1).

When $1\leq t\leq c(2-\alpha)^{-1}$ we use the inequality
$$N\bigl(\sqrt{q}B_2^n, tZ_q(\mu)\bigr)\ls N\bigl(\sqrt{q}B_2^n, c_1 (2-\alpha)^{-1} Z_q(\mu)\bigr)\,
N\bigl(Z_q(\mu), c_1^{-1}(2-\alpha)t Z_q(\mu)\bigr).$$ Noticing that
the latter covering number is less than
$\bigl(1+c(2-\alpha)^{-1}t^{-1}\bigr)^n$ completes the proof for the
first covering number in (\ref{eq:covering-2}) after elementary
calculations.

The bound for the second covering number in (\ref{eq:covering-2})
follows from the duality of entropy theorem. $\hfill \Box$

\medskip

Using Proposition~\ref{prop:projZq} instead of
Proposition~\ref{prop:projZq-small} we get the following:
\begin{proposition}\label{prop:covering-large}
Let $\mu $ be an isotropic log-concave measure on $\mathbb R^n$.
Assume that $2\ls q\ls n$. Then, for any $1\ls\alpha <2$ and any
$$1\ls t\ls\min\biggl\{\sqrt{q},\;
c_2 (2-\alpha)^{-1}L_n\left(\frac{n}{q}\right)^{\frac{\alpha+4}{2\alpha}}\biggr\}$$
we have
\begin{eqnarray}
&&\max\Bigl\{\log N\bigl(\sqrt{q}B_2^n,\; tZ_q(\mu )\bigr),\ \log
N\bigl(\sqrt{q}Z_q^\circ(\mu),\; tB_2^n\bigr)\Bigr\}\nonumber\\[1ex]
&&\hspace*{10em}\ls c(\alpha)L_n^{\frac{2\alpha}{\alpha+4}}\, \frac{n}{t^{\frac{2\alpha }{\alpha +4}}}
\,\max\left\{\log\frac{\sqrt{2q}}{t},\;\log\frac{L_n}{(2-\alpha)t}\right\},
\label{eq:covering-2-large}
\end{eqnarray}
where
$c(\alpha)\ls C(2-\alpha )^{-2/3}$
and $c_1, c_2, C>0$ are absolute constants.
\end{proposition}

\noindent \textit{Proof.}
We proceed along the same lines as before, and the only other thing that we need to take
into account is the fact that, for every $\varepsilon\in (0,1)$, $L_{\varepsilon n}\ls cL_n$ for some absolute constant $c$.
Note that here we are allowed to consider any $t$ up to $\sqrt{q}$
if we restrict ourselves to those $q$ that do not exceed $\sqrt{L_n}n^{3/4}$. $\hfill\Box$

\begin{remark}\label{rem:covering-large}\upshape
If we do not use the monotonicity of $L_n$ but we rather use the
bound $L_{\varepsilon n}\ls \sqrt[4]{\varepsilon n}$, we will end up
with an upper bound of the form
\begin{equation}\label{eq:covering-2-large-no-monotonicity}
 \frac{C}{(2-\alpha)^{\frac{4\alpha}{\alpha+8}}}
\frac{n^{\frac{2\alpha+8}{\alpha+8}}\log q}{t^{\frac{4\alpha
}{\alpha +8}}}
\end{equation}for $\max\{\log N\bigl(\sqrt{q}B_2^n, tZ_q(\mu )\bigr), \log
N\bigl(\sqrt{q}Z_q^\circ(\mu), tB_2^n\bigr) \}$.
We thus get a better exponent of $t$ for each $\alpha$, but the restrictions in the
proof force $t$ to be in the range
$$c(2-\alpha)^{-1}\sqrt[4]{n}\ls t\ls c(2-\alpha)^{-1}
\frac{n^{\frac{2\alpha+8}{4\alpha}}}{q^{\frac{\alpha+8}{4\alpha}}}$$
(again $t$ can vary up to $\sqrt{q}$ if, for example, $q\ls
n^{6/7}$). In this range, the latter bound for the covering numbers is
more efficient only if $L_n$ depends ``badly'' on $n$.
\end{remark}

\section{Upper bound for $M(Z_q(\mu ))$}

To make use of the covering estimates we have just obtained so as to give an
upper bound for $M(Z_q(\mu ))$, we employ the Dudley--Fernique decomposition
(see e.g.~\cite[\S 2.5.2]{Gian}). We consider the symmetric convex body
$K:=\sqrt{q}Z_q^{\circ }(\mu )$ and, for any $1\ls j\ls \log q$, we
consider the entropy number $N(K,2^{-j}R\,B_2^n)$, where $R=R(K)\ls
\sqrt{q}$. There exists $N_j\subseteq K$ with
$|N_j|=N(K,2^{-j}R\,B_2^n)$ such that, for any $x\in K$, there exists
$v\in N_j$ satisfying $\|x-v\|_2\ls 2^{-j}R$. We set $N_0=\{0\}$ and
$Z_j=N_j-N_{j-1}$. Then, we have:

\begin{lemma}\label{lem:5.1}
Let $K$ be a symmetric convex body in $\mathbb R^n$. For any $m\in
\mathbb N$ and any $x\in K$ there exist $z_1,\ldots,z_m,w_m$ with
$z_j\in Z_j\cap \frac{3R}{2^j}\,B_2^n$ and $w_m\in
\frac{R}{2^m}\,B_2^n$ such that
\begin{equation}\label{eq:5.1}x=z_1+\cdots +z_m+w_m,
\end{equation}where $R=R(K)$ is the radius of $K$.
\end{lemma}

\begin{theorem}[version for small $q$]\label{th:8.2}
Let $\mu $ be an isotropic log-concave measure on $\mathbb R^n$. For
every
$1\ls \alpha <2$,
for every
\begin{equation}\label{condition:5.3}
1\ls q\ls c (2-\alpha)^{-1/6} n^{\frac{\alpha+4}{3\alpha+8}},
\end{equation}
we have
\begin{equation}\label{eq:5.3}
M\bigl(Z_q(\mu )\bigr)\ls
C(2-\alpha)^{-1/3}
\frac{\sqrt{\max\bigl\{\log q,\;\log(2-\alpha)^{-1}\bigr\} } }{q^{\frac{\alpha}{2\alpha + 8}}},
\end{equation} where $C>0$ is an absolute constant. In particular, for every $1\ls q\ls n^{3/7}$,
\begin{equation}\label{eq:5.3a} M\bigl(Z_q(\mu )\bigr)
\ls C\frac{(\log q)^{5/6}}{\sqrt[6]{q}}.\end{equation}
\end{theorem}

\noindent \textit{Proof.} We set $K=\sqrt{q}Z_q^{\circ }(\mu )$. Using
Lemma \ref{lem:5.1}, for any $m\in \mathbb N$ and any $x\in K$ we
can find $(z_j)_{j\ls m}\subset Z_j\cap \frac{3R}{2^j}B_2^n$ and
$w_m\in \frac{R}{2^m}B_2^n$ such that $x=z_1+\cdots +z_m+w_m$. For
any $\theta\in S^{n-1}$ one has
\begin{equation}\label{eq:5.5} |\langle x,\theta\rangle|\ls \sum_{j=1}^m |\langle
z_j,\theta \rangle |+|\langle w_m,\theta\rangle|.
\end{equation} We write $\overline{z}=z/\| z\|_2$ for all $z\neq 0$.
We have
\begin{align}\label{eq:5.6}
w(K) & =\int_{S^{n-1}}\max_{x\in K}|\langle x,\theta
  \rangle|\,d\sigma(\theta) \\
  & \ls \sum_{j=1}^m\int_{S^{n-1}}\max_{z\in
  Z_j}|\langle \theta,z\rangle|\, d\sigma(\theta)+\int_{S^{n-1}}\max_{w\in 2^{-m}R\,B_2^n}|\langle
  w,\theta\rangle|\, d\sigma(\theta)\nonumber \\
  &\ls \sum_{j=1}^m\frac{3R}{2^j} \int_{S^{n-1}} \max_{z\in Z_j}|\langle \theta,\bar
  z\rangle|\, d\sigma(\theta) +\frac{R}{2^m} \nonumber \\
  &\ls \sum_{j=1}^m\frac{c_3R}{2^j}\frac{\sqrt{\log
  |Z_j|}}{\sqrt{n}} +\frac{R}{2^m} ,\nonumber
\end{align} where we have used the following:

\smallskip

\noindent \textit{Fact.} For any $u_1,\ldots, u_N \in S^{n-1}$ we have
\begin{equation}\label{eq:5.7}
\int_{S^{n-1}}\max_{j\ls N}|\langle \theta ,u_j\rangle|\,
d\sigma(\theta)\ls c_3\frac{\sqrt{\log N}}{\sqrt{n}}.
\end{equation}

\noindent By the definition of $Z_j$ and Proposition
\ref{prop:covering} we obtain
\begin{equation}\label{eq:5.8}\log|Z_j|\ls \log|N_j|+\log|N_{j-1}|\ls
c(\alpha )n \left (\frac{2^j}{R}\right )^{\frac{2\alpha
}{\alpha +4}} \max\bigl\{\log q,\;\log(2-\alpha)^{-1}\bigr\},
\end{equation}
where we assume that $R/2^m\geqslant 1$ and $c(\alpha )\ls C(2-\alpha )^{-\frac{2\alpha }{\alpha +4}}$.
Plugging this into (\ref{eq:5.6}) we conclude that
\begin{align}\label{eq:5.9}w(K) &\ls
\frac{c_5}{(2-\alpha )^{\frac{\alpha }{\alpha +4}}}\ R^{\frac{4}{\alpha +4}}\,
\sqrt{\max\biggl\{\log q,\;\log\frac1{2-\alpha}\biggr\}}
\ \sum_{j\ls m} \frac{1}{2^{\frac{4 j}{\alpha +4}}}+\frac{R}{2^m}\\
\nonumber &\ls \frac{C}{(2-\alpha )^{\frac{\alpha }{\alpha
+4}}}\ q^{\frac{2}{\alpha +4}}\,\sqrt{\max\biggl\{\log q,\;\log\frac1{2-\alpha}\biggr\}},
\end{align}
if $m$ is large enough so that
$R/2^m\simeq 1$.
It remains to observe that
\begin{equation*}
w(K)=\sqrt{q}w\bigl(Z_q^{\circ }(\mu)\bigr)=\sqrt{q}M\bigl(Z_q(\mu )\bigr),
\end{equation*}
and hence
\begin{equation*}
M\bigl(Z_q(\mu )\bigr)\ls \frac{C}{(2-\alpha )^{\frac{\alpha }{\alpha
+4}}}\frac{\sqrt{\max\bigl\{\log q,\;\log(2-\alpha)^{-1}\bigr\}}}{q^{\frac{\alpha}{2\alpha
+8}}}.
\end{equation*}
Finally, to obtain (\ref{eq:5.3a}), we set $\alpha =2-\frac{1}{\log q}$. $\hfill \Box$

\medskip

Given a symmetric isotropic convex body $K$ in ${\mathbb R}^n$
we apply Theorem \ref{th:8.2} with the isotropic measure $\mu_K$, that
has density $L_K^n\mathbf{1}_{{K/L_K}}$, and with $q= n^{3/7}$
(which is the optimal choice for this purpose), and we get:

\begin{theorem}\label{th:8.3}
Let $K\subset {\mathbb R}^n$ be isotropic and symmetric.
Then,
\begin{equation*}M(K)\ls C\frac{(\log n)^{5/6}}{L_K\sqrt[14]{n}}.\end{equation*}
\end{theorem}

Using Proposition \ref{prop:covering-large} instead of Proposition \ref{prop:covering}, we also get:

\begin{theorem}
Let $\mu $ be an isotropic log-concave measure on $\mathbb R^n$. For
every
$1\ls \alpha <2$,
for every
\begin{equation}\label{condition:th8.4.1}
1\ls q\ls c_1(2-\alpha)^{-1/3}\, L_n^{\frac{2\alpha}{2\alpha+4}}\, n^{\frac{\alpha+4}{2\alpha+4}},
\end{equation}
we have
\begin{equation}\label{eq:th8.4.1}
M\bigl(Z_q(\mu )\bigr)\ls C(2-\alpha)^{-1/3}\, L_n^{\frac{\alpha}{\alpha+4}}
\,\frac{\sqrt{\max\bigl\{\log q,\;\log L_n(2-\alpha)^{-1}\bigr\}}}{q^{\frac{\alpha}{2\alpha + 8}}},
\end{equation} where $c_1$, $C>0$ are absolute constants.
In particular, for every $q$ such that
$L_n\log q\ls q\ls \sqrt{L_n}\,n^{3/4}$,
\begin{equation}\label{eq:th8.4.2} M\bigl(Z_q(\mu )\bigr)\ls C\frac{\sqrt[3]{L_n}(\log q)^{5/6}}{\sqrt[6]{q}}\end{equation}
and, for every symmetric isotropic convex body $K$ in ${\mathbb R}^n$,
\begin{equation*}M(K) \ls C\frac{\sqrt[4]{L_n} (\log n)^{5/6}}{L_K\sqrt[8]{n}}.\end{equation*}
\end{theorem}

\begin{remark}\rm Similarly, using Remark \ref{rem:covering-large}
instead of Proposition \ref{prop:covering}, we see that, for every
$1\ls \alpha <2$,
for every
\begin{equation}\label{condition:th8.4.2}
c_1(2-\alpha)^{-2}\sqrt{n} \ls q\ls c_2(2-\alpha)^{-1/3} n^{\frac{2\alpha+8}{3\alpha+8}},
\end{equation}
we have
\begin{equation}\label{eq:th8.4.3}
M\bigl(Z_q(\mu )\bigr)\ls C (2-\alpha)^{-2/5}\,\sqrt[4]{n}^{\frac{2\alpha}{\alpha+8}}\,\frac{\sqrt{\log q}}{q^{\frac{\alpha}{\alpha + 8}}},
\end{equation} where $c_1$, $c_2$ and $C>0$ are absolute constants, and, for every $\sqrt{n}\log n\ls q\ls n^{6/7}$,
\begin{equation}\label{eq:th8.4.4} M\bigl(Z_q(\mu )\bigr)\ls C\frac{\sqrt[10]{n}(\log n)^{9/10}}{\sqrt[5]{q}}.\end{equation}
Therefore, for every isotropic symmetric convex body $K$ in ${\mathbb R}^n$,
\begin{equation*}M(K) \ls C\frac{(\log n)^{9/10}}{L_K \sqrt[14]{n}}.\end{equation*}
\end{remark}

\section{Further observations}

In this last Section we collect a number of additional observations
on the geometry of the centroid bodies $Z_q(\mu )$.

\medskip

\noindent \textbf{1. Inradius of projections.} We first provide lower
bounds for $R(Z_q(\mu )\cap F)$ and $R(Z_q^{\circ }(\mu )\cap F)$;
actually, they hold true for every $1\ls k<n$ and any $F\in
G_{n,k}$. By duality, these estimates (combined with
e.g.~Proposition~\ref{prop:projZq-small}) determine the inradius of
$P_F(Z_q^{\circ }(\mu ))$ and $P_F(Z_q(\mu ))$. Our starting point
is the next proposition, which can be essentially found in
\cite{Giannopoulos-Milman-1998}.

\begin{proposition}\label{prop:7-2}Let $A$ be a symmetric convex body in ${\mathbb R}^n$. Assume
that there exists $\gamma\gr 1$ such that
\begin{equation*}N(B_2^n,tA)\ls\exp
\Big(\frac{\gamma n}{t^p}\Big)\end{equation*}for every $t\gr 1$. For
every $\delta\in (0,1)$ and every $F\in G_{n,\lfloor \delta n\rfloor
}$ we have
\begin{equation*}w(A\cap F)\gr \frac{c\delta^{1/p}}{\gamma^{1/p}}.\end{equation*}
\end{proposition}

\noindent \textit{Proof.} Let $k=\lfloor \delta n\rfloor $ and consider
any $F\in G_{n,k}$. Using the assumption and the duality of entropy
theorem we see that the projection $P_F(A^{\circ })$ of $A^{\circ }$
onto $F$ satisfies
\begin{equation*}
N\bigl(P_F(A^{\circ }),tB_F\bigr)\ls N(A^{\circ },tB_2^n)\ls \exp
\left (\frac{\gamma k}{\delta t^p}\right ),\end{equation*}for every
$t\gr 1$. We apply Theorem \ref{th:section-alpha-regular} with
$W=P_F(A^{\circ })$, $n=k$ and $\varepsilon =1/2$. There exists
$H\in G_{k,\lfloor k/2\rfloor }(F)$ for which
\begin{equation*}P_F(A^{\circ })\cap H\subseteq \frac{c\gamma^{1/p}}{\delta^{1/p}}\, B_H.\end{equation*}
Taking polars in $H$ we see that $P_H(A\cap F) \supseteq
\frac{c\delta^{1/p}}{\gamma^{1/p}}B_H$.

Recall now that, given a symmetric convex body $C$ in ${\mathbb
R}^m$ and an $s$-dimensional subspace $L$ of ${\mathbb R}^m$, one
has $M(C\cap L) \lesssim \sqrt{m/s}\, M(C)$ (see \cite[Section
4.2]{Giannopoulos-Milman-handbook}). Therefore, setting $C = (A\cap
F)^{\circ }$ and $L = H$, we obtain
\begin{align*}w(A\cap F) &=M\bigl((A\cap F)^{\circ }\bigr)\gr\frac{1}{\sqrt{2}}
M\bigl((A\cap F)^{\circ }\cap H\bigr)=\frac{c}{\sqrt{2}}w\bigl(P_H(A\cap F )\bigr) \\
&\gr \frac{c'\delta^{1/p}}{\gamma^{1/p}},\end{align*}
as claimed. $\hfill\Box $

\begin{theorem}\label{th:8-2}Let $\mu $ be an isotropic log-concave measure on
${\mathbb R}^n$. Assume that $q\ls \sqrt{n}$. Then, for any
$1\ls\alpha <2$, $0<\delta <1$ and any $F\in G_{n,\lfloor\delta
n\rfloor }$ we have
\begin{equation*}
R\bigl(Z_q(\mu )\cap F\bigr)\gr w\bigl(Z_q(\mu )\cap F\bigr)
\gr c\frac{(2-\alpha )\delta^{\frac{\alpha +4}{2\alpha
}}}{\left(\log\frac1{2-\alpha}\right)^{\frac{\alpha+4}{2\alpha}}}
\,\sqrt{q}.\end{equation*}
\end{theorem}

\noindent \textit{Proof.} From Proposition \ref{prop:covering} we know
that
\begin{equation*}
\log N\bigl(\sqrt{q}B_2^n, tZ_q(\mu )\bigr)\ls
c(\alpha)\frac{n}{t^{\frac{2\alpha }{\alpha +4}}} \max\left\{\log
\frac{\sqrt{2q}}{t}, \;\log\frac1{(2-\alpha)t}\right\},
\end{equation*}
where $c(\alpha )\ls C(2-\alpha )^{-2/3}$. Thus, we may apply
Proposition \ref{prop:7-2} with
$$\gamma =c(\alpha)\frac{\max\{\log q, \log\frac1{2-\alpha}\}}{\sqrt{q}^p}$$
and $p=\frac{2\alpha }{\alpha +4}$. $\hfill\Box $

\medskip

A similar argument applies to $Z_q^{\circ }(\mu )$. Since
$M(Z_q^{\circ }(\mu ))=w(Z_q(\mu ))\simeq\sqrt{q}$ for all
$q\ls\sqrt{n}$, from the dual Sudakov inequality we have
\begin{equation*}
\log N\left(B_2^n, t\sqrt{q}Z_q^{\circ }(\mu )\right)\ls
\frac{cn}{t^2}
\end{equation*}for all $t\gr 1$. Applying Proposition \ref{prop:7-2} with
${\gamma=c q}$ and $p=2$ we get:

\begin{theorem}\label{th:8-3}Let $\mu $ be an isotropic log-concave measure on
${\mathbb R}^n$. Assume that $q\ls \sqrt{n}$. Then, for any
$1\ls\alpha <2$, $0<\delta <1$ and any $F\in G_{n,\lfloor\delta
n\rfloor }$ we have
\begin{equation*}
R\bigl(Z_q^{\circ }(\mu )\cap F\bigr)\gr
w\bigl(Z_q^{\circ }(\mu )\cap F\bigr)\gr \frac{c\sqrt{\delta }}{\sqrt{q}}.\end{equation*}
\end{theorem}

\smallskip

\noindent \mathversion{bold}\textbf{2. Upper bound for $M_{-k}(Z_q(\mu ))$.}\mathversion{normal}
Let $C$ be a
symmetric convex body in ${\mathbb R}^n$. For every $p\neq 0$, one
can define
\begin{equation}M_p(C):=\left(\int_{S^{n-1}}\|\theta\|_C^p d\sigma (\theta )\right)^{1/p}.\end{equation}
Litvak, Milman and Schechtman proved in \cite{LMS} that if $b$ is
the smallest constant for which $\| x\|\ls b\|x\|_2$ holds true for
every $x\in {\mathbb R}^n$, then
\begin{equation*}\max\bigg\{ M(C),c_1\frac{b\sqrt{q}}{\sqrt{n}}\bigg\}\ls M_p(C)\ls
\max\bigg\{ 2M(C),c_2\frac{b\sqrt{q}}{\sqrt{n}}\bigg\}
\end{equation*} for all $p\in [1,n]$, where $c_1,c_2>0$ are absolute
constants. In particular, \begin{equation}\label{eq:LMS}M_p(C)\simeq
M(C)\end{equation} as long as $p\ls k(C):=k_{\ast }(C^{\circ })$.
Klartag and Vershynin defined in \cite{Klartag-Vershynin} the
parameter
\begin{equation*}d(C)=\min\left\{ -\log \sigma\left (\left
\{ x\in S^{n-1}:\| x\|\ls\frac{M(C)}{2}\right\}\right
),n\right\}\end{equation*}and they observed that $d(C)$ is always
larger than $k(C)$. Their main result is an analogue of
(\ref{eq:LMS}) for negative values of $p$: one has
\begin{equation}\label{eq:KV}M_{-p}(C)\simeq M(C)\end{equation}as
long as $0<p\ls d(C)$.

Let $\mu $ be an isotropic log-concave measure on ${\mathbb R}^n$.
Since $M_{-p}(Z_q(\mu ))$ is clearly smaller than $M(Z_q(\mu ))$,
our next aim is to provide upper bounds for these quantities and to
compare them to the ones from Section 7. We will use a formula for
$M_{-k}(C)$ which appears in \cite{PaourisTAMS}.

\begin{lemma}\label{lem:formulaM-k}Let $C$ be a symmetric convex body in ${\mathbb R}^n$.
For every integer $1\ls k<n$,
\begin{equation}\label{W-kC}
M_{-k}(C) \simeq \left(\int_{G_{n,k}} {\rm vrad}\bigl(P_F(C^{\circ})\bigr)^{-k} d\nu_{n,k}(F) \right)^{-1/k}.
\end{equation}
\end{lemma}

The proof is simple: using the Blaschke-Santal\'{o} and the reverse
Santal\'{o} inequality (see~\cite{BM}) one can write
\begin{align}
M_{-k}^{-1}(C) &= \left( \int_{S^{n-1}}
\frac{1}{\| x\|_C^k} d\sigma(x) \right)^{1/k}\\
\nonumber &= \left(\frac{1}{\omega_{k}} \int_{G_{n,k}} \omega_{k}
\int_{S_{F}} \frac{1}{\|x\|_{C\cap F}^{k}}
d\sigma(x) d\nu_{n,k}(F)\right)^{1/k} \\
\nonumber &= \left(\int_{G_{n,k}} \frac{ |C\cap F|}{|B_{2}^{k}|}
d\nu_{n,k}(F) \right)^{1/k} \\
\nonumber &\simeq  \left(\int_{G_{n,k}}
\frac{|B_{2}^{k}|}{|P_F(C^{\circ })|} d\nu_{n,k}(F) \right)^{1/k}.
\end{align}

\begin{proposition}\label{prop:upper-bound-forM-kZq-small}Let $\mu $ be an
isotropic log-concave measure on ${\mathbb R}^n$. For every $2\ls
q\ls \sqrt{n}$ and $k\gr q^2\log^2q$ one has
\begin{equation}M_{-k}\bigl(Z_q(\mu )\bigr)\ls \frac{c\log^3q}{\sqrt{q}}\left (\frac{n}{k}\right
)^{3/2}.\end{equation}
\end{proposition}

\noindent \textit{Proof.} We choose $\alpha =2-\frac{1}{\log q}$. From
the proof of Proposition~\ref{prop:covering} we see that if
$t\leq\sqrt{q}$ and
\begin{equation*}\frac{q^2}{(\log q)^{2/3}}t^{2/3}\ls
n\end{equation*}then
\begin{equation*}
N\bigl(\sqrt{q}Z_q^\circ(\mu),\; tB_2^n\bigr)\ls \exp\left
(c_1\frac{n(\log q)^{2/3}}{t^{2/3}}\log\frac{2q}{t^2}\right ).
\end{equation*}
where $c_1$ is an absolute constant. Thus, for any integer $1\ls
k<n$ and any $F\in G_{n,k}$ we have
\begin{align*}
\sqrt{q}\Bigl|P_F\bigl(Z_q^{\circ }(\mu )\bigr)\Bigr|^{1/k} &\ls
t|B_2^k|^{1/k}\,N\Bigl(\sqrt{q}P_F\bigl(Z_q^{\circ }(\mu )\bigr), tB_F\Bigr)^{1/k}\\
&\ls t|B_2^k|^{1/k}\exp\left (c_1\frac{n(\log
q)^{2/3}}{kt^{2/3}}\log\frac{2q}{t^2}\right ).
\end{align*}Choosing $t\simeq (n/k)^{3/2}\log^3q$ we conclude that
\begin{equation*}\left (\frac{\Bigl|P_F\bigl(Z_q^{\circ }(\mu )\bigr)\Bigr|}{|B_2^k|}\right
)^{1/k}\ls \frac{c\log^3q}{\sqrt{q}}\left (\frac{n}{k}\right
)^{3/2},\end{equation*}and the result follows from Lemma
\ref{lem:formulaM-k}. $\hfill\Box $

\smallskip

\noindent \textit{Note.} Since $Z_q(\mu )\supseteq B_2^n$, one has the
obvious upper bound $M_{-k}(Z_q(\mu ))\ls 1$ for all $1\ls k<n$. So,
the estimate of Proposition \ref{prop:upper-bound-forM-kZq-small} is
non-trivial provided that $k\gr n\log^2q/q^{1/3}$.

\begin{remark}\rm Similarly, starting from Proposition
\ref{prop:covering-large} and following the proof of Proposition
\ref{prop:upper-bound-forM-kZq-small} we obtain the
following.\end{remark}

\begin{proposition}\label{prop:upper-bound-forM-kZq}Let $\mu $ be an
isotropic log-concave measure on ${\mathbb R}^n$. For every $2\ls
q\ls n$ and $k\gr L_n^{2/3}q$ one has
\begin{equation}M_{-k}\bigl(Z_q(\mu )\bigr)\ls \frac{cL_n}{\sqrt{q}}\left (\frac{n}{k}\right
)^{3/2}.\end{equation}
\end{proposition}

In view of (\ref{eq:KV}) a natural question is to give lower bounds
for $d(Z_q(\mu ))$; because, if $d(Z_q(\mu ))$ is large enough so
that we can use Proposition~\ref{prop:upper-bound-forM-kZq-small} or
Proposition~\ref{prop:upper-bound-forM-kZq}, then we would have an
alternative source of, possibly better, information on $M(Z_q(\mu
))$ as well. What we know is a simple lower bound for $k(Z_q(\mu ))$
and $k(Z_q^{\circ }(\mu ))$ in the range $2\ls q\ls q_{\ast}(\mu )$.
For every $2\ls q\ls q_{\ast}(\mu )$ one has
\begin{equation}\min\Bigl\{ k\bigl(Z_q(\mu )\bigr),k\bigl(Z_q^{\circ
}(\mu )\bigr)\Bigr\}\gr \frac{c_1n}{q},\end{equation}where $c_1$ is an
absolute constant.

To see this, recall that $B_2^n\subseteq Z_q(\mu )\subseteq
c_2qB_2^n$, and hence
\begin{equation*}R\bigl(Z_q(\mu )\bigr)\ls c_2q\quad\hbox{and}\quad R\bigl(Z_q^{\circ
}(\mu )\bigr)\ls 1.\end{equation*}Then, using (\ref{eq:wZq-small}) we get
\begin{equation*}k_{\ast }\bigl(Z_q(\mu )\bigr)\gr c_3
n\frac{w^2\bigl(Z_q(\mu )\bigr)}{R^2\bigl(Z_q(\mu )\bigr)}\gr
c_4n\frac{q}{q^2}=\frac{c_4n}{q}.\end{equation*} Also, using the
fact that $w(C^{\circ})w(C)\gr 1$ for every symmetric convex body
$C$, and taking into account (\ref{eq:wZq-small}), we see that
$w(Z_q^{\circ }(\mu ))\gr c/\sqrt{q}$ for every $2\ls q\ls
q_{\ast}(\mu )$, and thus
\begin{equation*}k_{\ast }\bigl(Z_q^{\circ }(\mu )\bigr)\gr c_5
n\frac{w^2\bigl(Z_q^{\circ }(\mu )\bigr)}{R^2\bigl(Z_q^{\circ }(\mu )\bigr)}\gr
\frac{c_6n}{q}.\end{equation*}

\medskip

\noindent \textbf{3. Small ball probability estimates.} Given a
symmetric convex body $C$ in ${\mathbb R}^n$ we set
$\overline{C}=|C|^{-1/n}C$. In this last subsection we describe an
approach which can lead to small ball probability estimates for the
centroid bodies. It is convenient to normalize the volume, and
consider $\overline{Z_q}(\mu )$ instead of $Z_q(\mu )$. Recall that
if $q\ls \sqrt{n}$ then $|Z_q(\mu )|^{1/n}\simeq \sqrt{q/n}$, and
hence
\begin{equation*}\overline{Z_q}(\mu )\simeq\sqrt{n/q}\,Z_q(\mu
).\end{equation*} Then, $w(\overline{Z_q}(\mu
))\simeq\sqrt{n/q}w(Z_q(\mu ))\simeq\sqrt{n}$, which implies
\begin{equation*}
\log N\bigl(\overline{Z_q}(\mu ),sB_2^n\bigr)\ls c_1n\left
(\frac{w\bigl(\overline{Z_q}(\mu )\bigr)}{s}\right
)^2\ls\frac{c_2n^2}{s^2}.\end{equation*} We use the following fact
(\cite{DP}, Lemma~5.6).

\begin{lemma}\label{lem:key-lemma-dp}
Let $C$ be a centered convex body of volume $1$ in ${\mathbb R}^n$.
Assume that, for some $s>0$,
\begin{equation}\label{eq:assumption-key-lemma-dp}r_{s}:=\log{N(K,sB_2^n)} < n.\end{equation} Then,
\begin{equation*}I_{-r_s}(K) \ls c s.\end{equation*}
\end{lemma}
\noindent We apply Lemma \ref{lem:key-lemma-dp} for $\overline{Z_q}(\mu )$ to
get:

\begin{proposition}\label{prop:I-qZq}Let $\mu $ be an isotropic
log-concave measure on ${\mathbb R}^n$. If $2\ls q\ls\sqrt{n}$ then
\begin{equation*}I_{-r}\bigl(\overline{Z_q}(\mu )\bigr)\ls
\frac{c_3n}{\sqrt{r}}\end{equation*}for all
$1\ls r\ls cn$.
\end{proposition}

From Markov's inequality we have the small ball probability estimate
\begin{equation*}\biggl|\Bigl\{ x\in \overline{Z_q}(\mu ):\| x\|_2\ls
\varepsilon I_{-r}\bigl(\overline{Z_q}(\mu )\bigr)\Bigr\}\biggr|\ls
\varepsilon^r.\end{equation*}

\bigskip

\medskip

\footnotesize

\bigskip

\bigskip

\footnotesize

\noindent \textsc{Apostolos Giannopoulos:} Department of Mathematics,
University of Athens, Panepistimioupolis 157\,84, Athens, Greece.

\noindent \textit{E-mail:} \texttt{apgiannop@math.uoa.gr}

\medskip

\noindent \textsc{Pantelis Stavrakakis:} Department of Mathematics,
University of Athens, Panepistimioupolis 157\,84, Athens, Greece.

\noindent \textit{E-mail:} \texttt{pantstav@yahoo.gr}

\medskip

\noindent \textsc{Antonis Tsolomitis:} Department of Mathematics,
University of the Aegean, Karlovassi 832\,00, Samos, Greece.

\noindent \textit{E-mail:} \texttt{antonis.tsolomitis@gmail.com}

\medskip

\noindent \textsc{Beatrice-Helen Vritsiou:} Department of Mathematics,
University of Athens, Panepistimioupolis 157\,84, Athens, Greece.

\noindent \textit{E-mail:} \texttt{bevritsi@math.uoa.gr}

\medskip


\begin{thebibliography}{100}
\footnotesize

\bibitem{Artstein-Milman-Szarek-2004}
\textrm{S.\ Artstein, V.\ D.\ Milman and S.\ J.\ Szarek},
\textsl{Duality of metric entropy},
Annals of Math.\ \textbf{159} (2004), 1313--1328.

\bibitem{Ball}
\textrm{K.\ M.\ Ball},
\textsl{Logarithmically concave functions and sections of convex sets in ${\mathbb R}^n$},
Studia Math.\ \textbf{88} (1988), 69--84.

\bibitem{Bastero-2007}
\textrm{J.\ Bastero},
\textsl{Upper bounds for the volume and diameter of $m$-dimensional sections of convex bodies},
Proc.\ Amer.\ Math.\ Soc.\ \textbf{135} (2007), 1851--1859.

\bibitem{Bou}
\textrm{J.\ Bourgain},
\textsl{On the distribution of polynomials on high dimensional convex sets},
Geom.\ Aspects of Funct.\ Analysis (Lindenstrauss-Milman eds.), Lecture Notes in
Math.\ \textbf{1469} (1991), 127--137.

\bibitem{BM}
\textrm{J.\ Bourgain and V.\ D.\ Milman},
\textsl{New volume ratio properties for convex symmetric bodies in ${\mathbb R}^n$},
Invent.\ Math.\ \textbf{88} (1987), 319--340.

\bibitem{BLM}
\textrm{J.\ Bourgain, J.\ Lindenstrauss and V.\ D.\ Milman},
\textsl{Minkowski sums and symmetrizations},
Geom.\ Aspects of Funct.\ Analysis (Lindenstrauss-Milman eds.), Lecture Notes in
Math.~\textbf{1317} (1988), 44--74.

\bibitem{BGVV-book-isotropic}
\textrm{S.\ Brazitikos, A.\ Giannopoulos, P.\ Valettas and B-H.\ Vritsiou},
\textsl{Notes on isotropic convex bodies},
Book in preparation, available at {\tt http://users.uoa.gr/\~{}apgiannop/}.

\bibitem{DP}
\textrm{N.\ Dafnis and G.\ Paouris},
\textsl{Small ball probability estimates, $\psi_2$-behavior and the hyperplane conjecture},
J.\ Funct.\ Anal.\ \textbf{258} (2010), 1933--1964.


\bibitem{Gian}
\textrm{A.\ Giannopoulos},
\textsl{Notes on isotropic convex bodies},
Lecture Notes, Warsaw 2003, available at {\tt http://users.uoa.gr/\~{}apgiannop/}.

\bibitem{Giannopoulos-Milman-1997}
\textrm{A.\ Giannopoulos and V.\ D.\ Milman},
\textsl{On the diameter of proportional sections of a symmetric convex body},
International Mathematics Research Notices (1997) No.\ 1, 5--19.

\bibitem{Giannopoulos-Milman-1998}
\textrm{A.\ A.\ Giannopoulos and V.\ D.\ Milman},
\textsl{Mean width and diameter of proportional sections of a symmetric convex body},
J.\ Reine Angew.\ Math.\ \textbf{497} (1998), 113--139.

\bibitem{Giannopoulos-Milman-handbook}
\textrm{A.\ Giannopoulos and V.\ D.\ Milman},
\textsl{Euclidean structure in finite-dimensional normed spaces},
Handbook of the Geometry of Banach Spaces (Johnson-Lindenstrauss eds.), Vol.\ 1 (2001), 707--779.

\bibitem{Giannopoulos-Milman-Tsolomitis-2005}
\textrm{A.\ Giannopoulos, V.\ D.\ Milman and A.\ Tsolomitis},
\textsl{Asymptotic formulas for the diameter of sections of symmetric convex bodies},
Journal of Functional Analysis \textbf{223} (2005), 86--108.

\bibitem{Giannopoulos-Paouris-Valettas-2011}
\textrm{A.\ Giannopoulos, G.\ Paouris and P.\ Valettas},
\textsl{On the existence of subgaussian directions for log-concave measures},
Contemporary Mathematics \textbf{545} (2011), 103--122.

\bibitem{Giannopoulos-Paouris-Valettas-2012}
\textrm{A.\ Giannopoulos, G.\ Paouris and P.\ Valettas},
\textsl{On the distribution of the $\psi_2$-norm of linear functionals on
isotropic convex bodies}, Geom.\ Aspects of Funct.\ Analysis,
Lecture Notes in Math. \textbf{2050} (2012), 227--253.

\bibitem{GPV-JFA}
\textrm{A.\ Giannopoulos, G.\ Paouris and B-H.\ Vritsiou},
\textsl{A remark on the slicing problem},
Journal of Functional Analysis 262 (2012), 1062--1086.

\bibitem{Gordon}
\textrm{Y.\ Gordon},
\textsl{On Milman's inequality and random subspaces which escape through a mesh in $\textbf{R}^n$},
Lecture Notes in Mathematics \textbf{1317} (1988), 84--106.

\bibitem{Kl}
\textrm{B.\ Klartag},
\textsl{On convex perturbations with a bounded isotropic constant},
Geom.\ Funct.\ Anal.~\textbf{16} (2006), 1274--1290.

\bibitem{KM1}
\textrm{B.\ Klartag and E.\ Milman},
\textsl{Centroid bodies and the logarithmic Laplace transform---A unified approach},
J.\ Funct.\ Anal.\ \textbf{262} (2012), 10--34.

\bibitem{KM2}
\textrm{B.\ Klartag and E.\ Milman},
\textsl{Inner regularization of
log-concave measures and small-ball estimates}, Geom.\ Aspects of
Funct.\ Analysis, Lecture Notes in Math. \textbf{2050} (2012),
267--278.

\bibitem{KlM1}
\textrm{B.\ Klartag and V.\ D.\ Milman},
\textsl{Rapid Steiner symmetrization of most of a convex body and the slicing problem},
Combin.\ Probab.\ Comput.\ \textbf{14}, no. 5--6 (2005) 829--843.

\bibitem{Klartag-Vershynin}
B.\ Klartag and R.\ Vershynin,
\textsl{Small ball probability and Dvoretzky theorem},
Israel J.\ Math.~\textbf{157} (2007), 193--207.

\bibitem{LMS}
\textrm{A.\ Litvak, V.\ D.\ Milman and G.\ Schechtman},
\textsl{Averages of norms and quasi-norms},
Math.\ Ann.\ \textbf{312} (1998), 95--124.

\bibitem{Litvak-Milman-Pajor-Tomczak-2006}
\textrm{A.\ Litvak, V.\ D.\ Milman, A.\ Pajor and N.\ Tomczak-Jeagermann},
\textsl{Entropy extension},
Funct.\ Anal.\ Appl.\ \textbf{40} (2006), 298--303.

\bibitem{LYZ}
\textrm{E.\ Lutwak, D.\ Yang and G.\ Zhang},
\textsl{$L^p$ affine isoperimetric inequalities},
J.\ Differential Geom.~\textbf{56} (2000), 111--132.

\bibitem{Milman-1988}
\textrm{V.\ D.\ Milman},
\textsl{Isomorphic symmetrization and geometric inequalities},
Lecture Notes in Mathematics \textbf{1317} (1988), 107--131.

\bibitem{Milman-1991}
\textrm{V.\ D.\ Milman},
\textsl{Some applications of duality relations},
Lecture Notes in Mathematics \textbf{1469} (1991), 13--40.

\bibitem{Milman-Pajor-1989}
\textrm{V.\ D.\ Milman and A.\ Pajor},
\textsl{Isotropic position and inertia ellipsoids and zonoids of the
  unit ball of a normed $n$-dimensional space},
Lecture Notes in Mathematics \textbf{1376}, Springer, Berlin (1989), 64--104.

\bibitem{Milman-Schechtman-book}
\textrm{V.\ D.\ Milman and G.\ Schechtman},
\textsl{Asymptotic theory of finite dimensional normed spaces},
Lecture Notes in Mathematics \textbf{1200} (1986), Springer, Berlin.

\bibitem{PaourisGAFA}
\textrm{G.\ Paouris},
\textsl{Concentration of mass in convex bodies},
Geometric and Functional Analysis  \textbf{16} (2006), 1021--1049.

\bibitem{PaourisTAMS}
\textrm{G.\ Paouris},
\textsl{Small ball probability estimates for
log-concave measures}, Trans.\ Amer.\ Math.\ Soc.~\textbf{364} (2012),
287--308.

\bibitem{Pisier-1989}
\textrm{G.\ Pisier},
\textsl{A new approach to several results of V.\ Milman},
J.\ Reine Angew.\ Math.\ \textbf{393} (1989), 115--131.

\bibitem{Pisier-book}
\textrm{G.\ Pisier},
\textsl{The volume of convex bodies and Banach space geometry},
Cambridge Tracts in Mathematics \textbf{94} (1989).

\bibitem{Schneider-book}
\textrm{R.\ Schneider},
\textsl{Convex bodies: the Brunn-Minkowski theory},
Encyclopedia of Mathematics and its Applications \textbf{44},
Cambridge University Press, Cambridge (1993).

\bibitem{Valettas-private}
\textrm{P.\ Valettas},
\textsl{Upper bound for the $\ell $-norm in the isotropic position},
Private communication.

\bibitem{Vershynin-2006}
\textrm{R.\ Vershynin},
\textsl{Isoperimetry of waists and local versus global asymptotic convex geometries}
(with an appendix by M. Rudelson and R. Vershynin),
Duke Mathematical Journal \textbf{131} (2006), 1--16.

\bibitem{Vritsiou-IMRN}
\textrm{B-H.\ Vritsiou},
\textsl{Further unifying two approaches to the hyperplane conjecture},
Int.\ Math.\ Res.\ Not.\ (2012), \texttt{DOI: 10.1093/imrn/rns263}.

\end{thebibliography}
\end{document}